\documentclass[twosided,reqno,11pt]{amsart}

\usepackage{amsmath,latexsym,amssymb,pstcol}
\usepackage[ps,dvips,all]{xy}
\SelectTips{cm}{} 

\def\sqr#1#2{{\vcenter{\hrule height.#2pt
        \hbox{\vrule width.#2pt height#1pt \kern#1pt
                \vrule width.#2pt}
        \hrule height.#2pt}}}
\def\square{\mathchoice\sqr64\sqr64\sqr{4}3\sqr{3}3}
\def\QED{\hfill$\square$}
\def\demo{\noindent{\bf Proof: }}
\def\split{\genfrac{}{}{0pt}1}

\newtheorem{Theorem}{\bf Theorem}[section]
\newtheorem{Lemma}[Theorem]{\bf Lemma}
\newtheorem{Corollary}[Theorem]{\bf Corollary}
\newtheorem{Proposition}[Theorem]{\bf Proposition}

\newtheorem{Remark}[Theorem]{\bf Remark}
\newtheorem{Example}[Theorem]{\bf Example}
\newtheorem{Definition}[Theorem]{\bf Definition}

\newcommand {\ZZ}{\mathbb{Z}}
\newcommand {\NN}{\mathbb{N}}
\newcommand {\PP}{\mathbb{P}}

\newcommand{\mif}{\mbox{if} ~}
\newcommand{\s}{\; | \;}

\DeclareMathOperator{\reg}{reg}

\DeclareMathOperator{\Proj}{Proj}

\begin{document}

\baselineskip=16.5pt

\title[ Monomial and toric ideals associated to Ferrers graphs]
{\Large\bf Monomial and toric ideals associated to \\ Ferrers graphs}

\author[A. Corso and U. Nagel]
{Alberto Corso \and Uwe Nagel}

\thanks{AMS 2000 {\em Mathematics Subject Classification}.
Primary:  05A15, 13D02, 13D40, 14M25; Secondary: 05C75, 13C40, 13H10,  14M12,
 52B05.}


\address{Department of Mathematics, University of Kentucky, Lexington,
Kentucky 40506}
\email{corso@ms.uky.edu}
\email{uwenagel@ms.uky.edu}

\vspace{.1in}

\begin{abstract}
Each partition $\lambda = (\lambda_1, \lambda_2, \ldots, \lambda_n)$
determines a so-called Ferrers tableau or, equivalently, a Ferrers
bipartite graph. Its edge ideal, dubbed Ferrers ideal, is a
squarefree monomial ideal that is generated by quadrics. We show
that such an ideal has a $2$-linear minimal free resolution, i.e.\
it defines a small subscheme. In fact, we prove that this property
characterizes Ferrers graphs among bipartite graphs. Furthermore,
using a method of Bayer and Sturmfels, we provide an explicit
description of the maps in its minimal free resolution: This is
obtained by associating a suitable polyhedral cell complex to the
ideal/graph. Along the way, we also determine the irredundant
primary decomposition of any Ferrers ideal. We conclude our analysis
by studying several features of toric rings of Ferrers graphs. In
particular we recover/establish formul\ae\/ for the Hilbert series,
the Castelnuovo-Mumford regularity, and the multiplicity of these rings.
While most of the previous works in this highly investigated area of
research involve path counting arguments, we offer here a new and
self-contained approach based on results from Gorenstein liaison
theory.
\end{abstract}

\maketitle

\vspace{0.5in}

\section{Introduction}

\noindent A {\it Ferrers graph} is a bipartite graph on two distinct
vertex sets ${\bf X}=\{ x_1, \ldots, x_n\}$ and ${\bf Y} =\{ y_1,
\ldots, y_m \}$ such that if $(x_i, y_j)$ is an edge of $G$, then so
is $(x_p, y_q)$ for $1 \leq p \leq i$ and $1 \leq q \leq j$. In
addition, $(x_1,y_m)$ and $(x_n,y_1)$ are required to be edges of
$G$. For any Ferrers graph $G$ there is an associated sequence of
non-negative integers $\lambda = (\lambda_1, \lambda_2, \ldots,
\lambda_n)$, where $\lambda_i$ is the degree of the vertex $x_i$.
Notice that the defining properties of a Ferrers graph imply that
$\lambda_1=m \geq \lambda_2 \geq \cdots \geq \lambda_n \geq 1$; thus
$\lambda$ is a {\it partition}. Alternatively, we can associate to a
Ferrers graph a diagram ${\mathbf T}_{\lambda}$, dubbed {\it Ferrers
tableau}, consisting of an array of $n$ rows of cells with
$\lambda_i$ adjacent cells, left justified, in the $i$-th row.

Ferrers graphs/tableaux have a prominent place in the literature as
they have been studied in relation to chromatic polynomials
\cite{Brenti,EvW}, Schubert varieties \cite{Ding,Develin},
hypergeometric series \cite{Haglund}, permutation statistics
\cite{Butler,EvW}, quantum mechanical operators \cite{Varvak},
inverse rook problems \cite{GJW,Ding,Develin,Mitchell}. More
generally, algebraic and combinatorial aspects of bipartite graphs
have been studied in depth (see, e.g., \cite{SVV, HH} and the
comprehensive monograph  \cite{vila}). In this paper, which is the
first of a series \cite{CN2, CN}, we are interested in the algebraic
properties of the {\it edge ideal} $I=I(G)$ and the {\it toric ring}
$K[G]$ associated to a Ferrers graph $G$. The edge ideal is the
monomial ideal of the polynomial ring $R = K[x_1, \ldots, x_n, y_1,
\ldots, y_m]$ over the field $K$ that is generated by the monomials
of the form $x_iy_j$, whenever the pair $(x_i,y_j)$ is an edge of
$G$. $K[G]$ is instead the monomial subalgebra generated by the
elements $x_iy_j$. An example is illustrated in Figure 1:
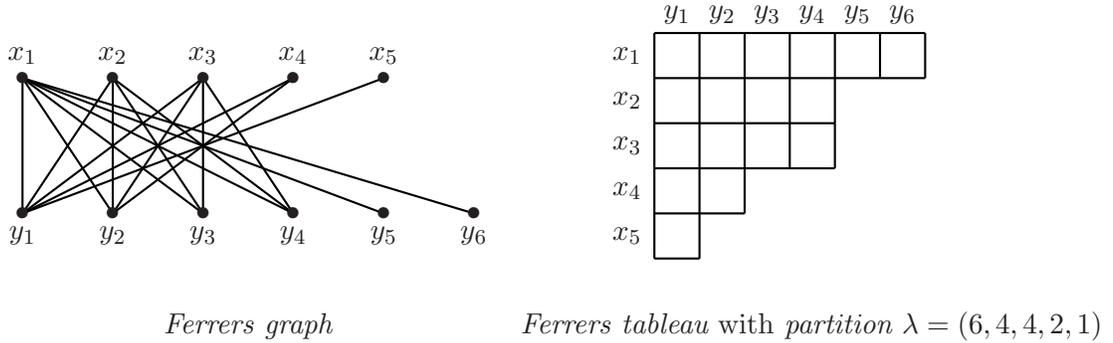
\begin{center}
\begin{pspicture}(-3,-1.7)(20,2.5)
\psset{xunit=.6cm, yunit=.6cm}

\psline[linestyle=solid](0,3)(0,0)
\psline[linestyle=solid](0,3)(2,0)
\psline[linestyle=solid](0,3)(4,0)
\psline[linestyle=solid](0,3)(6,0)
\psline[linestyle=solid](0,3)(8,0)
\psline[linestyle=solid](0,3)(10,0)
\psline[linestyle=solid](2,3)(0,0)
\psline[linestyle=solid](2,3)(2,0)
\psline[linestyle=solid](2,3)(4,0)
\psline[linestyle=solid](2,3)(6,0)
\psline[linestyle=solid](4,3)(0,0)
\psline[linestyle=solid](4,3)(2,0)
\psline[linestyle=solid](4,3)(4,0)
\psline[linestyle=solid](4,3)(6,0)
\psline[linestyle=solid](6,3)(0,0)
\psline[linestyle=solid](6,3)(2,0)
\psline[linestyle=solid](8,3)(0,0)

\rput(0,3.5){$x_1$} \rput(2,3.5){$x_2$} \rput(4,3.5){$x_3$}
\rput(6,3.5){$x_4$} \rput(8,3.5){$x_5$}

\rput(0,-.5){$y_1$} \rput(2,-.5){$y_2$} \rput(4,-.5){$y_3$}
\rput(6,-.5){$y_4$} \rput(8,-.5){$y_5$} \rput(10,-.5){$y_6$}

\rput(0,3){$\bullet$} \rput(2,3){$\bullet$} \rput(4,3){$\bullet$}
\rput(6,3){$\bullet$}\rput(8,3){$\bullet$}

\rput(0,0){$\bullet$} \rput(2,0){$\bullet$} \rput(4,0){$\bullet$}
\rput(6,0){$\bullet$} \rput(8,0){$\bullet$} \rput(10,0){$\bullet$}

\rput(5,-2.5){{\it Ferrers graph}}

\psline[linestyle=solid](14,4)(20,4)
\psline[linestyle=solid](14,3)(20,3)
\psline[linestyle=solid](14,2)(18,2)
\psline[linestyle=solid](14,1)(18,1)
\psline[linestyle=solid](14,0)(16,0)
\psline[linestyle=solid](14,-1)(15,-1)
\psline[linestyle=solid](14,4)(14,-1)
\psline[linestyle=solid](15,4)(15,-1)
\psline[linestyle=solid](16,4)(16,0)
\psline[linestyle=solid](17,4)(17,1)
\psline[linestyle=solid](18,4)(18,1)
\psline[linestyle=solid](19,4)(19,3)
\psline[linestyle=solid](20,4)(20,3)

\rput(13.4,3.5){$x_1$} \rput(13.4,2.5){$x_2$} \rput(13.4,1.5){$x_3$}
\rput(13.4,0.5){$x_4$} \rput(13.4,-0.5){$x_5$}

\rput(14.5,4.4){$y_1$} \rput(15.5,4.4){$y_2$} \rput(16.5,4.4){$y_3$}
\rput(17.5,4.4){$y_4$} \rput(18.5,4.4){$y_5$} \rput(19.5,4.4){$y_6$}

\rput(17.5,-2.5){{\it Ferrers tableau} with {\it partition}
$\lambda=(6,4,4,2,1)$}
\end{pspicture}
\end{center}
\[
I= (x_1y_1, x_1y_2, x_1y_3, x_1y_4, x_1y_5, x_1y_6, x_2y_1, x_2y_2,
x_2y_3, x_2y_4, x_3y_1, x_3y_2, x_3y_3, x_3y_4, x_4y_1, x_4y_2,
x_5y_1)
\]
\begin{center}
{\bf Figure 1: Ferrers graph, tableau and ideal}
\end{center}

Throughout this article $\lambda=(\lambda_1, \ldots, \lambda_n)$
will always denote a fixed partition associated to a Ferrers graph
$G_\lambda$ with corresponding Ferrers ideal $I_{\lambda}$. In
Section 2 we describe several fine numerical invariants attached to
the ideal $I_{\lambda}$. In Theorem~\ref{betti} we show that each
Ferrers ideal defines a small subscheme in the sense of Eisenbud,
Green, Hulek, and Popescu \cite{EGHP}, i.e.\ the free resolution of
$I_{\lambda}$ is $2$-linear. More precisely, we give an explicit
--- but at the same time surprisingly simple --- formula for the
Betti numbers of the ideal $I_{\lambda}$; namely, we show that:
\[
\beta_i(R/I_{\lambda}) = {\lambda_1 \choose i} + {\lambda_2+1
\choose i}+ {\lambda_3+2 \choose i} + \ldots + {\lambda_n+n-1
\choose i} - {n \choose i+1}
\]
for $1 \leq i \leq \max\{ \lambda_i+i-1 \}$. Furthermore, the Hilbert
series is:
\[
\sum_{k \geq 0} \dim_K [R/I_{\lambda}]_k \cdot t^k = \frac{1}{(1-t)^m} +
\frac{t}{(1-t)^{m+n+1}}
\cdot \sum_{j=1}^n (1-t)^{\lambda_j + j}.
\]
Notice that
the formula for the Betti numbers involves a minus sign: This is
quite an unusual
phenomenon for Betti numbers, as they tend, in general, to have an
enumerative interpretation. In order to determine the Betti numbers
it is essential to find a $($not necessarily irredundant$)$ primary
decomposition of $I_{\lambda}$. We refine this decomposition into an
irredundant one in Corollary~\ref{decomposition}, where we observe,
in particular, that the number of prime components is related to the
outer corners of the Ferrers tableau. For instance, in the case of
the ideal $I_{\lambda}$ described in Figure 1 we have that it is the
intersection of $5$ $(=4$ outer corners $+1)$ components:
\[
I_{\lambda}=(y_1, \ldots , y_6) \cap (x_1, y_1, y_2, y_3, y_4) \cap
(x_1, x_2, x_3, y_1, y_2) \cap (x_1, x_2, x_3, x_4, y_1) \cap (x_1,
\ldots, x_5).
\]
We conclude Section 2 by identifying, in terms of the shape of the
tableau, the unmixed $($Corollary~\ref{unmixed}$)$ and
Cohen-Macaulay $($Corollary~\ref{CM}$)$ members in the family of
Ferrers ideals. The latter result also follows from recent work of
Herzog and Hibi \cite{HH}.

There are relatively few general classes of ideals for which an
explicit minimal free resolution is known: The most noteworthy such
families include the Koszul complex, the Eagon-Northcott complex
\cite{EN}, and the resolution of generic monomial ideals \cite{BS}
$($see also \cite{BPS}$)$. In Section 3 we analyze even further the
minimal free resolution of a Ferrers ideal $I_{\lambda}$ and obtain
a surprisingly elegant description of the differentials in the
resolution in Theorem~\ref{exactness}.  In some sense, this is a
prototypical result as it provides the minimal free resolution of
several classes of ideals obtained from Ferrers ideals  by
appropriate specializations of the variables $($see \cite{CN2} for
further details$)$. Our description of the free resolution of a
Ferrers ideal relies on the theory of cellular resolutions as
developed by Bayer and Sturmfels in \cite{BS} $($see also
\cite{MS}$)$. More precisely, let $\Delta_{n-1} \times \Delta_{m-1}$
denote the product of two simplices of dimensions $n-1$ and $m-1$,
respectively. Given a Ferrers ideal $I_{\lambda}$, we associate to
it the polyhedral cell complex $X_{\lambda}$ consisting of the faces
of $\Delta_{n-1} \times \Delta_{m-1}$ whose vertices are labeled by
generators of $I_{\lambda}$ $($see Definition~\ref{def-cell-c}$)$.
By the theory of Bayer and Sturmfels, $X_{\lambda}$ determines a
complex of free modules. Using an inductive argument we show in
Theorem~\ref{exactness} that this complex is in fact the multigraded
minimal free resolution of the ideal $I_{\lambda}$. While leaving
the details to the main body of the paper, we illustrate the
situation in the case of the partition $\lambda=(4,3,2,1)$, which is
the largest we can draw. In this case the polyhedral cell complex
$X_{\lambda}$ can actually be identified with the subdivision of the
simplex $\Delta_3$ pictured below $($see \cite{CN2} for additional
details$)$:
\begin{center}
\begin{pspicture}(-8,-2.5)(0,3.5)
\psset{xunit=.6cm, yunit=.6cm}

\psline[linestyle=solid](-7,3)(-3,3)
\psline[linestyle=solid](-7,2)(-3,2)
\psline[linestyle=solid](-7,1)(-4,1)
\psline[linestyle=solid](-7,0)(-5,0)
\psline[linestyle=solid](-7,-1)(-6,-1)
\psline[linestyle=solid](-7,3)(-7,-1)
\psline[linestyle=solid](-6,3)(-6,-1)
\psline[linestyle=solid](-5,3)(-5,0)
\psline[linestyle=solid](-4,3)(-4,1)
\psline[linestyle=solid](-3,3)(-3,2)

\rput(-7.4,2.5){$x_1$} \rput(-7.4,1.5){$x_2$} \rput(-7.4,0.5){$x_3$}
\rput(-7.4,-0.5){$x_4$}

\rput(-6.5,3.4){$y_1$} \rput(-5.5,3.4){$y_2$} \rput(-4.5,3.4){$y_3$}
\rput(-3.5,3.4){$y_4$}

\psline[linewidth=1mm, linestyle=solid](-7,3)(-7,0)
\psline[linewidth=1mm, linestyle=solid](-6,0)(-6,1)
\psline[linewidth=1mm, linestyle=solid](-5,1)(-5,2)
\psline[linewidth=1mm, linestyle=solid](-4,2)(-4,3)

\psline[linewidth=1mm, linestyle=solid](-7.08,3)(-3.92,3)
\psline[linewidth=1mm, linestyle=solid](-5.08,2)(-3.92,2)
\psline[linewidth=1mm, linestyle=solid](-6.08,1)(-4.92,1)
\psline[linewidth=1mm, linestyle=solid](-7.08,0)(-5.92,0)

\rput(-2,1){$\leadsto$}
\end{pspicture}
\begin{pspicture}(0,-2.5)(7,3.5)
\psset{xunit=.8cm, yunit=.8cm}

\pscustom[linestyle=dashed,fillstyle=solid,fillcolor=lightgray]{
\psline(2,2)(4.75,.5) \psline[liftpen=0](4.75,-1.5)(2,0) }

\psline[linestyle=dashed]{*-*}(0,0)(4,0)
\psline[linestyle=dashed](4,0)(4,4)

\pscustom[linestyle=dashed,fillstyle=solid,fillcolor=gray]{
\psline{*-*}(2,2)(2.75,-1.5) \psline[liftpen=0](2.75,-1.5)(2,0) }

\pscustom[linestyle=dashed,fillstyle=solid,fillcolor=gray]{
\psline{*-*}(2,2)(4,2) \psline[liftpen=0](4,2)(4.75,.5) }

\psline[linestyle=dashed]{*-*}(2,2)(2,0)
\psline[linestyle=dashed](4,0)(5.5,-3)
\psline[linestyle=solid](0,0)(5.5,-3)
\psline[linestyle=solid](0,0)(4,4)
\psline[linestyle=solid]{*-*}(4,4)(5.5,-3)

\psline[linestyle=solid](2,2)(2.75,-1.5)
\psline[linestyle=solid](2,2)(4.75,0.5)

\rput(4.75,-1.5){$\bullet$} \rput(4.75,.5){$\bullet$}

\rput(1.6,2.3){$x_1y_1$} \rput(1.6,-0.3){$x_2y_1$}
\rput(2.35,-1.85){$x_3y_1$} \rput(-0.55,0){$x_4y_1$}

\rput(4,4.3){$x_1y_4$} \rput(3.5,2.3){$x_1y_3$}
\rput(5.4,0.7){$x_1y_2$} \rput(3.5,0.3){$x_2y_3$}

\rput(5.5,-3.3){$x_3y_2$} \rput(4.25,-1.75){$x_2y_2$}

\end{pspicture}
\end{center}
\begin{center}
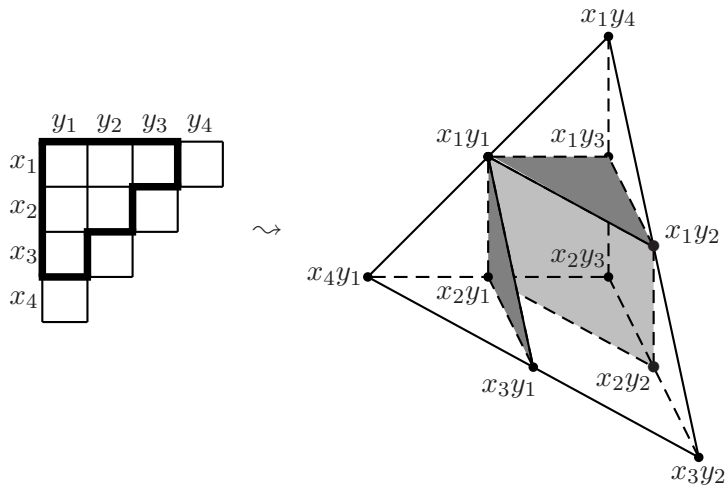
{\bf Figure 2: Ferrers tableau and associated polyhedral cell
complex}
\end{center}
In particular, we observe that $X_{\lambda}$ has four
3-dimensional cells: Two of them are isomorphic to $\Delta_3$
whereas the remaining two are isomorphic to either $\Delta_1
\times \Delta_2$ or $\Delta_2 \times \Delta_1$. A grey shading in
the picture above also indicates how the polyhedral cell complex
corresponding to the partition $(3,2,1)$ sits inside
$X_{\lambda}$.

In Section 4 we prove the converse of Theorem~\ref{betti}. Namely,
we show that any edge ideal of a bipartite graph with a $2$-linear
resolution necessarily arises from a Ferrers graph $($see
Theorem~\ref{converse}$)$. One of the ingredients of the proof is a
well-known characterization of edge ideals of graphs with a
$2$-linear resolution in terms of complementary graphs, due to Fr\"oberg
 \cite{Froberg} (see also \cite{EGHP2}).

The starting point of Section~\ref{Toric} is the observation that
the {\it toric ring} of a Ferrers graph can be identified with a
special {\it ladder determinantal ring}. We then proceed to
recover/establish formul\ae\/ for the Hilbert series and other
invariants associated with these rings. We remark that this is a
highly investigated part of mathematics that has been the subject of
the work of  many researchers. Among the extensive, impressive and
relevant literature we single out \cite{Abhy, BH, Cn, Cn2, CnH, G,
HT, K, KM, KP, KR, Kulkarni, N, R, wang}. While most of these works
involve --- to a different extent --- path counting arguments, we
offer here a new and self-contained approach that yields easy proofs
of explicit formul\ae\/ for the Hilbert series, the
Castelnuovo-Mumford regularity, and the multiplicity of the toric
rings of Ferrers graphs. This method, which is based on results from
Gorenstein liaison theory (see \cite{migbook} for a comprehensive
introduction), has been pioneered in \cite{KMMNP}, where
it was proved that every standard determinantal ideal is {\it
glicci}, i.e.\ it is in the Gorenstein liaison class of a complete
intersection (see also \cite{MNR}). Recently, Gorla \cite{Gorla} has
considerably refined these arguments to show that all ladder
determinantal ideals are glicci. This result can be used to
establish first a simple recursive formula, which we then turn into
an explicit formula that involves only positive summands.

\section{Betti numbers and primary decompositions of Ferrers ideals}

The main result of this section, Theorem~\ref{betti}, provides a
$($not necessarily irredundant$)$ primary decomposition of a Ferrers
ideal $I_{\lambda}$ as well as the Betti numbers of its minimal free
resolution. A particularly relevant situation is given by a maximal
bipartite graph, in which case the Ferrers tableau has a rectangular
shape of size $n \times m$ and $I_{\lambda} = (x_1, \ldots, x_n) (y_1,
\ldots, y_m)$; using a variety of techniques including determinantal
ideals, residual intersections and Gr\"obner basis, \cite{CrVV} and
\cite{BG} describe additional features of this and other related
subideals in connection with the so-called Dedekind-Mertens Lemma. A
hook-shaped tableau is the other extremal case; in this situation
$I_{\lambda} = x_1(y_1, \ldots, y_m) + y_1(x_1, \ldots, x_n)$.

For sake of simplicity we will  denote the partition associated to a
Ferrers graph by $\lambda = (\lambda_1, \lambda_2, \ldots,$
$\lambda_s, 1, \ldots, 1)$ with $\lambda_1 \geq \lambda_2 \geq
\cdots \geq \lambda_s \geq 2$. Furthermore, we denote the {\it dual
partition} by $\lambda^* = (\lambda_1^*, \lambda_2^*, \ldots,
\lambda_m^*)$, where $\lambda_j^*$ is the degree of the vertex
$y_j$. Observe that $\lambda_1 = m$, $\lambda_1^*=n$, and $1 \leq s
= \lambda_2^* \leq n$. We also recall that the Hilbert series of the
graded $K$-algebra $R/I_{\lambda}$ is:
\[
P(R/I_{\lambda}, t) := \sum_{k \geq 0} h_{R/I_{\lambda}}(k) \cdot
t^k := \sum_{k \geq 0} \dim_K [R/I_{\lambda}]_k \cdot t^k,
\]
where $h_{R/I_{\lambda}}$ is the Hilbert function of
$R/I_{\lambda}$. It is well-known that this series is a
rational function.

\begin{Theorem}\label{betti}
Let $G$ be a Ferrers graph with associated partition $\lambda =
(\lambda_1, \lambda_2, \ldots, \lambda_s,$ $1, \ldots, 1)$ and let
$I_{\lambda}$ be the edge ideal in $K[x_1, \ldots, x_n, y_1, \ldots,
y_m]$ associated to $G$. Then a $($not necessarily irredundant$)$
primary decomposition of $I_{\lambda}$ is$:$
\[
(y_1, \ldots, y_{\lambda_1}) \cap (x_1, y_1, \ldots,
y_{\lambda_2}) \cap (x_1, x_2, y_1, \ldots, y_{\lambda_3}) \cap
\ldots \cap \hspace{2cm}
\]
\[
\hspace{3cm} \cap (x_1, \ldots, x_{s-1}, y_1, \ldots,
y_{\lambda_s}) \cap (x_1, \ldots, x_s, y_1) \cap (x_1, \ldots,
x_n)
\]
and the minimal $\ZZ$-graded free resolution of $_{\lambda}$  is
$2$-linear with $i${\rm{-th}} Betti number given by$:$
\[
\beta_i(R/I_{\lambda}) = {\lambda_1 \choose i} + {\lambda_2+1
\choose i}+ {\lambda_3+2 \choose i} + \ldots + {\lambda_n+n-1
\choose i} - {n \choose i+1}
\]
for $1 \leq i \leq \max_j \{ \lambda_j+j-1 \}$. Furthermore, the
Hilbert series is$:$
\[
P(R/I_{\lambda}, t) = \frac{1}{(1-t)^m} + \frac{t}{(1-t)^{m+n+1}}
\cdot \sum_{j=1}^n (1-t)^{\lambda_j + j}.
\]
\end{Theorem}

\demo We proceed by induction on $n$. If $n=1$ then
$\lambda=(\lambda_1)=(m)$, $s=1$ and $I_{\lambda} = x_1(y_1, \ldots,
y_m)=(x_1) \cap (y_1, \ldots, y_m) = (y_1, \ldots, y_m) \cap (x_1,
y_1) \cap (x_1)$. Moreover, the resolution of $I$ is given by $($a
shifted$)$ Koszul complex on $m$ generators. Hence the $i$-th Betti
number is
\[
\beta_i(R/I_{\lambda}) = {m \choose i} = {m \choose i} - {1 \choose
i+1},
\]
as the latter term is zero. Furthermore, the Hilbert series is
\[
P(R/I_{\lambda}, t) = \frac{1}{(1-t)^m} + \frac{t}{(1-t)},
\]
as claimed.

Suppose now $n \geq 2$. We distinguish two main cases:
$\lambda_n=1$ and $\lambda_n \geq 2$.

\noindent
We first deal with the case $\lambda_n=1$.
In addition to the partition $\lambda$ we
consider the partition $\lambda' = (\lambda_1, \ldots, \lambda_{n-1})$.
Notice that the index $s$ is the same for both $\lambda$
and $\lambda'$. By induction hypothesis we have that a primary
decomposition of $I_{\lambda'}$ is
\[
(y_1, \ldots, y_{\lambda_1}) \cap (x_1, y_1, \ldots,
y_{\lambda_2}) \cap (x_1, x_2, y_1, \ldots, y_{\lambda_3})
\cap \ldots \cap \hspace{2cm}
\]
\[
\hspace{3cm}
\cap (x_1, \ldots, x_{s-1}, y_1, \ldots, y_{\lambda_s})
\cap (x_1, \ldots, x_s, y_1) \cap (x_1, \ldots, x_{n-1}).
\]
Let $J$ denote the intersection of all the components in the above
primary decomposition that contain $x_ny_1$. Thus,
$I_{\lambda'} = J \cap (x_1, \ldots, x_{n-1})$. Now observe that
$I_{\lambda} = I_{\lambda'}+ (x_ny_1)$, thus using the above
primary decomposition for $I_{\lambda'}$, we get
\begin{eqnarray*}
I_{\lambda} & = &
J \cap (x_1, \ldots, x_{n-1}, x_ny_1) \\
& =& J \cap (x_1, \ldots, x_{n-1}, y_1) \cap (x_1,
\ldots, x_{n-1}, x_n) \\
&=& J \cap (x_1, \ldots, x_n),
\end{eqnarray*}
which is, after some inspection,
exactly the asserted primary decomposition of $I_{\lambda}$.

We now turn to the  Betti numbers of $I_{\lambda}$.  From the
given primary decomposition, it follows that $I_{\lambda'}
\colon x_ny_1 = (x_1, \ldots, x_{n-1})$. Hence we have the following
short exact sequence.
\[
0 \rightarrow R/(x_1, \ldots, x_{n-1})[-2]
\stackrel{\cdot x_ny_1}{\longrightarrow} R/I_{\lambda'}
\longrightarrow R/I_{\lambda} \rightarrow 0.
\]
Using a mapping cone construction we obtain that
\[
\beta_i(R/I_{\lambda}) = \beta_i(R/I_{\lambda'}) +
\beta_{i-1}(R/(x_1, \ldots, x_{n-1})) = \beta_i(R/I_{\lambda'}) +
{n-1 \choose i-1}.
\]
Hence, using our inductive assumption we have that
$\beta_i(R/I_{\lambda})$ is given by
\[
{\lambda_1 \choose i} + {\lambda_2+1 \choose i}
+  {\lambda_3+2 \choose i} + \ldots + {\lambda_{n-1}+n-2 \choose i}
- {n -1 \choose i+1} + {n-1 \choose i-1}.
\]
However, one can easily check the identity ${n-1 \choose i-1} - {n-1
\choose i+1} = {n \choose i} - {n \choose i+1}$, which provides the
expected form of $\beta_i(R/I_{\lambda})$.

Moreover, the above exact sequence provides for the Hilbert series
\[
P(R/I_{\lambda}, t) = P(R/I_{\lambda'}, t) - \frac{t^2}{(1-t)^{m+1}}.
\]
Using the induction hypothesis an easy computation provides the
claim for the Hilbert series of $R/I_{\lambda}$.

Suppose now $\lambda_n \geq 2$ and consider the
partition $\lambda'= (\lambda_1, \lambda_2, \ldots, \lambda_{n-1},
\lambda_n-1)$.
If in addition we assume that $\lambda_n \geq 3$, then the
index $s$ of the partition $\lambda'$ also equals $n$.
The inductive hypothesis provides that a primary
decomposition of $I_{\lambda'}$ is
\[
(y_1, \ldots, y_{\lambda_1}) \cap (x_1, y_1, \ldots,
y_{\lambda_2}) \cap (x_1, x_2, y_1, \ldots, y_{\lambda_3})
\cap \ldots \cap \hspace{2cm}
\]
\[
\hspace{3cm}
\cap (x_1, \ldots, x_{n-1}, y_1, \ldots, y_{\lambda_n-1})
\cap (x_1, \ldots, x_n, y_1) \cap (x_1, \ldots, x_n).
\]
Let $J$ denote the intersection of all the components in the above
primary decomposition that do not contain $x_ny_{\lambda_n}$.
Thus, $I_{\lambda'} = J \cap (x_1, \ldots, x_{n-1}, y_1, \ldots,
y_{\lambda_n-1})$. Now observe that $I_{\lambda} = I_{\lambda'}+
(x_n y_{\lambda_n})$, thus using the above primary decomposition for
$I_{\lambda'}$, we get
\begin{eqnarray*}
I_{\lambda} & = & J \cap (x_1, \ldots, x_{n-1},
y_1, \ldots, y_{\lambda_n-1}, x_ny_{\lambda_n}) \\
& =& J \cap (x_1, \ldots, x_{n-1}, y_1, \ldots, y_{\lambda_n-1},
y_{\lambda_n}) \cap (x_1,
\ldots, x_{n-1}, y_1, \ldots, y_{\lambda_n-1}, x_n) \\
&=& J \cap (x_1, \ldots, x_{n-1}, y_1, \ldots, y_{\lambda_n}),
\end{eqnarray*}
which is, after some inspection, exactly the asserted primary
decomposition of $I_{\lambda}$.
Turning to the  Betti numbers of $I_{\lambda}$,  the
given primary decomposition implies that $I_{\lambda'}
\colon x_ny_{\lambda_n} = (x_1, \ldots, x_{n-1}, y_1, \ldots,
y_{\lambda_n-1})$. Hence, by a similar mapping cone argument
as above, we obtain that
\[
\beta_i(R/I_{\lambda}) = \beta_i(R/I_{\lambda'}) +
\beta_{i-1}(R/(x_1, \ldots, x_{n-1}, y_1, \ldots, y_{\lambda_n-1}))
= \beta_i(R/I_{\lambda'}) + {\lambda_n-1+n-1 \choose i-1}.
\]
Therefore, using our inductive assumption we get that
$\beta_i(R/I_{\lambda})$ is given by
\[
{\lambda_1 \choose i} + \ldots + {\lambda_{n-1}+n-2 \choose i}
+ {\lambda_n-1+n-1 \choose i} -{n \choose i+1} + {\lambda_n-1+n-1
\choose i-1},
\]
which can be rewritten in the form
\[
{\lambda_1 \choose i} + \ldots + {\lambda_{n-1}+n-2 \choose i}
+ {\lambda_n+n-1 \choose i} -{n \choose i+1},
\]
as claimed.

To finish our proof, let us assume that $\lambda_n=2$ and consider
the partition $\lambda'=(\lambda_1, \ldots, \lambda_{n-1}, 1)$, whose
index $s$ is $n-1$. Moreover, by inductive assumption we have that
a primary decomposition of $I_{\lambda'}$ is
\[
(y_1, \ldots, y_{\lambda_1}) \cap (x_1, y_1, \ldots,
y_{\lambda_2}) \cap (x_1, x_2, y_1, \ldots, y_{\lambda_3})
\cap \ldots \cap \hspace{2cm}
\]
\[
\hspace{3cm}
\cap (x_1, \ldots, x_{n-2}, y_1, \ldots, y_{\lambda_{n-1}})
\cap (x_1, \ldots, x_{n-1}, y_1) \cap (x_1, \ldots, x_n).
\]
Let $J$ denote the intersection of all the components in the above
primary decomposition that  contain $x_ny_2$. Thus, $I_{\lambda'} =
J \cap (x_1, \ldots, x_{n-1},y_1)$. Now observe that $I_{\lambda} =
I_{\lambda'}+ (x_ny_2)$, thus using the above primary decomposition
for $I_{\lambda'}$, we get
\begin{eqnarray*}
I_{\lambda} & = & J \cap (x_1, \ldots, x_{n-1},
y_1, x_ny_2) \\
& =& J \cap (x_1, \ldots, x_{n-1}, y_1, y_2) \cap (x_1,
\ldots, x_{n-1}, y_1, x_n) \\
&=& J \cap (x_1, \ldots, x_{n-1}, y_1, y_2),
\end{eqnarray*}
which is, after some inspection, exactly the asserted primary
decomposition of $I_{\lambda}$.
The
given primary decomposition of $I_{\lambda'}$ provides that $I_{\lambda'}
\colon x_ny_2 = (x_1, \ldots, x_{n-1}, y_1)$.
Hence,  a mapping cone argument implies that
\[
\beta_i(R/I_{\lambda}) = \beta_i(R/I_{\lambda'}) +
\beta_{i-1}(R/(x_1, \ldots, x_{n-1}, y_1)) = \beta_i(R/I_{\lambda'})
+ {n \choose i-1}.
\]
Therefore, using our inductive assumption we see that
$\beta_i(I_{\lambda})$ is given by
\[
{\lambda_1 \choose i} + \ldots + {\lambda_{n-1}+n-2 \choose i}
+ {n \choose i} -{n \choose i+1} + {n
\choose i-1},
\]
which can be rewritten in the form
\[
{\lambda_1 \choose i} + \ldots + {\lambda_{n-1}+n-2 \choose i}
+ {n+1 \choose i} -{n \choose i+1},
\]
which gives the asserted formula as ${n+1 \choose i} = {2+(n-1) \choose i}$.

The claim about the Hilbert series follows similarly as in the case
$\lambda_n = 1$. We omit the details. \QED

\medskip

Theorem \ref{betti} allows us to compute further invariants of
Ferrers ideals:

\begin{Corollary} \label{cor-inv}
Adopt the notation of {\rm Theorem~\ref{betti}}. Then the height of
the edge ideal $I_{\lambda}$ of a Ferrers graph $G$ is $\min \{
\min_j \{ \lambda_j+j-1 \}, n \}$, the projective dimension of the
factor ring $R/I_{\lambda}$ is $\max_j\{ \lambda_j+j-1 \}$ and the
Castelnuovo-Mumford regularity ${\rm reg}(I_{\lambda})$ is equal to
$2$.
\end{Corollary}

\begin{Remark}
{\rm By a result of Herzog-Hibi-Zheng \cite{HHZ}, all the powers
$I_{\lambda}^k$ have a linear resolution so that the
Castelnuovo-Mumford regularity of $I_{\lambda}^k$, for any integer
$k \geq 1$, is ${\rm reg}(I_{\lambda}^k)=2k$. }
\end{Remark}

\begin{Example}\label{maximal-bipartite} {\rm
The edge ideal $I_{\lambda}$ of the complete bipartite graph  on $n$
and $m$ vertices, respectively, is $I_{\lambda} = (x_1, \ldots, x_n)
(y_1, \ldots, y_m) = (x_1, \ldots, x_n) \cap (y_1, \ldots, y_m)$ and
has $i${\rm{-th}} Betti number$:$
\[
\beta_i(R/I_{\lambda}) = {m+n \choose i+1} - {m \choose i+1} - {n
\choose i+1}
\]
for $1 \leq i \leq n+m-1$.
In fact, the formula for the Betti numbers follows from
Theorem~\ref{betti} and \cite[Formula (1.48)]{Gould}, or simply by
using the Mayer-Vietoris sequence. }
\end{Example}

\medskip

The primary decomposition of the Ferrers ideal $I_{\lambda}$
described in Theorem~\ref{betti} can be refined into an irredundant
one by using the shape of the Ferrers tableau ${\mathbf
T}_{\lambda}$. More precisely, define recursively indices $j_0,
\ldots, j_t$ by setting $j_0=0$ and, for $i\geq 0$,
\[
j_{i+1} = \max \{ k \,|\, \lambda_k = \lambda_{j_i +1 } \}.
\]
Note that $\lambda_{j_1}=m$, $j_t=n$, and that the pairs $(j_i,
\lambda_{j_i})$, $i=1, \ldots, t$, are the coordinates of the {\it
outer corners} of the Ferrers tableau ${\mathbf T}_{\lambda}$. In
addition, set $\lambda_{j_{t+1}}=0$ and, accordingly, $(x_1, \ldots,
x_{j_0}) = (0) = (y_1, \ldots, y_{\lambda_{j_{t+1}}})$. With this
notation, we state next our refinement of the primary decomposition
described in Theorem~\ref{betti}:

\begin{Corollary}\label{decomposition}
The irredundant primary decomposition of the Ferrers ideal
$I_{\lambda}$ is$:$
\[
I_{\lambda} = \bigcap_{i=1}^{t+1} (x_1, \ldots, x_{j_{i-1}}, y_1,
\ldots, y_{\lambda_{j_i}}),
\]
where the pairs $(j_i, \lambda_{j_i})$, $i=1, \ldots, t$, correspond
to the $t$ outer corners of the Ferrers tableau of $I_{\lambda}$. In
particular, $I_{\lambda}$ is the intersection of $t+1$ prime ideals.
\end{Corollary}
\demo This follows by inspecting the decomposition given in
Theorem~\ref{betti}. \QED

\begin{Corollary}\label{unmixed}
Adopt the notation of\ {\rm Theorem~\ref{betti}} as well as the one
established above. Then the following conditions are equivalent$:$
\begin{itemize}
\item[$({\it a})$]
$I_{\lambda}$ is unmixed\/$;$

\item[$({\it b})$]
$n=m$ and the inside corners $(j_{i-1}, \lambda_{j_i})$, for $i=2,
\ldots, t$, of the Ferrers tableau ${\mathbf T}_{\lambda}$ of
$I_{\lambda}$ lie on the main anti-diagonal of ${\mathbf
T}_{\lambda}$, i.e.\ on  $\{ (p,q) \,|\, p+q = m \}$.
\end{itemize}
\end{Corollary}
\demo The equivalence of the conditions follows immediately from
Corollary~\ref{decomposition}. \QED

\bigskip

Ferrers ideals are rarely Cohen-Macaulay. In fact, we get:

\begin{Corollary}\label{CM}  The following conditions are
equivalent$:$
\begin{itemize}
\item[$({\it a})$]
$I_{\lambda}$ is unmixed and connected in codimension one\/$;$

\item[$({\it b})$]
$n=m$ and $\lambda = (n, n-1, n-2, \ldots, 3, 2, 1);$

\item[$({\it c})$]
$I_{\lambda}$ is a Cohen-Macaulay ideal.
\end{itemize}
In particular, in this case the $i${\rm{-th}} Betti number is$:$
\[
\beta_i(R/I_{\lambda}) = i {n+1 \choose i+1}
\]
for $1 \leq i \leq n(=m)$. Moreover, the Cohen-Macaulay type is $n$
and the Hilbert series $P(R/I_{\lambda}, t)$ is given by$:$
\[
P(R/I_{\lambda},t) = \frac{1+nt}{(1-t)^n}.
\]
\end{Corollary}

\demo Corollary \ref{cor-inv} shows that $({\it b})$ implies $({\it
c})$. Condition $({\it a})$ is always a consequence of $({\it c})$.
Using Corollary~\ref{decomposition}, we see that $({\it a})$ implies
$({\it b})$.

Concerning the Betti numbers of $I_{\lambda}$, Theorem~\ref{betti}
and the shape of the partition $\lambda$ provide that
\[
\beta_i(R/I_{\lambda}) = n {n \choose i} - {n \choose i+1}
\]
for $1 \leq i \leq n$. On the other hand, an easy calculation shows
that the latter expression equals $i {n+1 \choose i+1}$, as claimed. The
statement about the Hilbert series follows immediately from
Theorem~\ref{betti}. \QED

\medskip

We note that the equivalence of conditions $({\it b})$ and $({\it
c})$ above can also be deduced from a recent result of  Herzog and
Hibi \cite[Theorem 3.4]{HH}. In general, the condition that a
projective subscheme $Z \subset \PP^n$ is equidimensional and
connected in codimension one is only a necessary condition for $Z$
being arithmetically Cohen-Macaulay. However, if $Z$ is defined by a
monomial ideal, then it seems often the case that this condition is
also sufficient.

\section{Minimal free resolutions of Ferrers ideals}

In this section we explicitly describe the minimal free resolution
of every Ferrers ideal. For our construction we use cellular
resolutions and polyhedral cell complexes as introduced  by Bayer
and Sturmfels in \cite{BS}. First, we briefly recall some basic
notions but we refer to \cite{BS} $($or \cite{MS}$)$ for a more
detailed introduction to the topic. A {\it polyhedral cell complex}
$X$ is a finite collection of convex polytopes $($in some ${\mathbb
R}^N)$ called faces $($or cells$)$ of $X$ such that:
\begin{enumerate}
\item
if $P \in X$ and $F$ is a face of $P$, then $F \in X$;

\item
if $P, Q \in X$ then $P \cap Q$ is a face of both $P$ and $Q$.
\end{enumerate}
Let $F_k(X)$ be the set of $k$-dimensional faces. Each cell complex
admits an {\it incidence function} $\varepsilon$ on $X$ where
$\varepsilon(Q,P) \in \{ 1, -1 \}$ if $Q$ is a facet of $P \in X$.
$X$ is called a {\it labeled} cell complex if each vertex $i$ has a
vector ${\mathbf a}_i \in {\mathbb N}^N$ $($or the monomial
${\mathbf z}^{{\mathbf a}_i}$, where ${\mathbf z}^{{\mathbf a}_i}$
denotes a monomial in the variables $z_1, \ldots, z_N)$ as label. The
label of an
arbitrary face $Q$ of $X$ is the exponent ${\mathbf a}_Q$, where
${\mathbf z}^{{\mathbf a}_Q} := {\rm lcm}\, ({\mathbf z}^{{\mathbf
a}_i} \,|\, i \in Q)$. Each labeled cell complex determines a
complex of free $R$-modules, where $R = K[z_1,\ldots,z_N]$. The
cellular complex ${\mathcal F}_X$ supported on $X$ is the complex of
free $\ZZ^N$-graded $R$-modules
\[
{\mathcal F}_X\colon \quad 0 \rightarrow S^{F_d(X)}
\stackrel{\partial_d}{\longrightarrow} S^{F_{d-1}(X)}
\stackrel{\partial_{d-1}}{\longrightarrow} \cdots
\stackrel{\partial_2}{\longrightarrow} S^{F_1(X)}
\stackrel{\partial_1}{\longrightarrow} S^{F_0(X)}
\stackrel{\partial_0}{\longrightarrow} S \rightarrow 0,
\]
where $d= \dim X$ and $S^{F_k(X)}:= \displaystyle \bigoplus_{P \in
F_k(X)} R[-{\mathbf a}_P]$. The map $\partial_k$ is defined by
\[
\partial_k(e_P) := \sum_{{Q {\rm \ facet \ of \ } P}}
\varepsilon(P,Q) \cdot {\mathbf z}^{{\mathbf a}_P -{\mathbf a}_Q}
\cdot e_Q,
\]
where $\{ e_P \,|\, P \in F_k(X) \}$ is a basis of $S^{F_k(X)}$ and
$e_{\emptyset} :=1 $. If ${\mathcal F}_X$ is acyclic, then it
provides a free $\ZZ^N$-graded resolution of the image $I$ of
$\partial_0$, that is the ideal generated by the labels of the
vertices of $X$. In  this case, ${\mathcal F}_X$ is called a {\it
cellular resolution} of $I$.

\medskip

We are ready to describe a cellular minimal free resolution for each
Ferrers ideal. First, let us consider the complete bipartite graph
${\mathcal K}_{n,m}$ that corresponds to the edge ideal $(x_1,
\ldots, x_n)(y_1, \ldots, y_m)$ $($or the partition $\lambda$, where
$\lambda_i=m$ for $i=1, \ldots, n)$. To this graph, we associate the
polyhedral cell complex $X_{n, m}$ given by the face complex of the
polytope $\Delta_{n-1} \times \Delta_{m-1}$ obtained by taking the
cartesian product of the $(n-1)$-simplex $\Delta_{n-1}$ and the
$(m-1)$-simplex $\Delta_{m-1}$. Labeling the vertices of
$\Delta_{n-1}$ by $x_1,\ldots,x_n$ and the ones of $\Delta_{m-1}$ by
$y_1,\ldots,y_m$, the vertices of the cell complex $X_{n,m}$ are
naturally labeled by the monomials $x_iy_j$ with $1 \leq i \leq n$
and $1 \leq j \leq m$. An easy example is illustrated below:
\begin{center}
\begin{pspicture}(0,0)(11,4)
\psset{xunit=.6cm, yunit=.6cm}

\rput(1.5,3.5){$x_1$} \rput(1.5,2.5){$x_2$} \rput(2.6,4.5){$y_1$}
\rput(3.6,4.5){$y_2$} \rput(4.6,4.5){$y_3$}

\rput(6.5,3){$\leadsto$}

\psline[linestyle=solid](2,4)(5,4)
\psline[linestyle=solid](2,3)(5,3)
\psline[linestyle=solid](2,2)(5,2)
\psline[linestyle=solid](2,4)(2,2)
\psline[linestyle=solid](3,4)(3,2)
\psline[linestyle=solid](4,4)(4,2)
\psline[linestyle=solid](5,4)(5,2)

\rput(9,2){$\bullet$} \rput(12,1){$\bullet$}
\rput(10.5,4.5){$\bullet$} \rput(14,3){$\bullet$}
\rput(17,2){$\bullet$} \rput(15.5,5.5){$\bullet$}

\rput(8.2,2.3){$x_1y_1$} \rput(11.9,0.5){$x_1y_2$}
\rput(9.7,4.8){$x_1y_3$} \rput(13.2,3.3){$x_2y_1$}
\rput(16.9,1.5){$x_2y_2$} \rput(14.7,5.8){$x_2y_3$}

\psline[linestyle=solid](9,2)(12,1)
\psline[linestyle=solid](9,2)(10.5,4.5)
\psline[linestyle=solid](12,1)(10.5,4.5)
\psline[linestyle=solid](17,2)(15.5,5.5)
\psline[linestyle=solid](12,1)(17,2)
\psline[linestyle=solid](10.5,4.5)(15.5,5.5)
\psline[linestyle=dashed](14,3)(17,2)
\psline[linestyle=dashed](14,3)(15.5,5.5)
\psline[linestyle=dashed](9,2)(14,3)
\end{pspicture}
\end{center}
\begin{center}
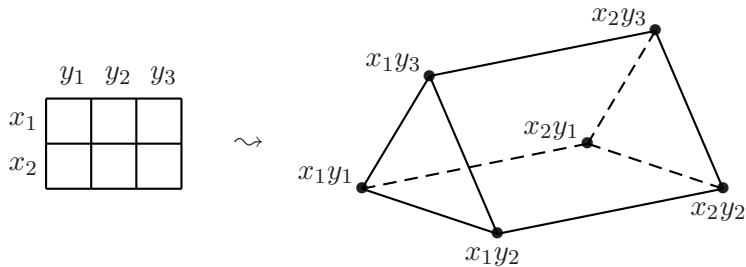
{\bf Figure 3: Topological viewpoint, $\Delta_1 \times \Delta_2$}
\end{center}
To simplify notation, we denote the monomial that labels the face $P
\in X_{n,m}$ by $m_P$. In general, we observe that $\deg m_P = \dim
P + 2$ for each face $P \in X_{n,m}$.

\medskip

We are now in the position to define the cell complex that will
support the cellular resolution of a given Ferrers ideal. As above,
we fix  a partition $\lambda=(\lambda_1, \lambda_2, \ldots,
\lambda_n)$ with $\lambda_1=m$, corresponding to a Ferrers graph
$G_{\lambda}$ and a Ferrers ideal $I _{\lambda} \subset R =
K[x_1,\ldots,x_n, y_1,\ldots,y_m]$.

\begin{Definition} \label{def-cell-c}
{\rm The polyhedral cell complex $X_{\lambda}$ associated to the
partition ${\lambda}$ is the labeled subcomplex of $X_{n,m}$ consisting of
 the faces of $X_{n,m}$ whose vertices are labeled by
all the monomials generating the Ferrers ideal $I_{\lambda}$. }
\end{Definition}

\noindent Using the Ferrers tableau ${\mathbf T}_{\lambda}$ we get a
more explicit, yet simple, description of the cell complex
$X_{\lambda}$. In fact, it is easy to see that the  facets
of $X_{\lambda}$ are in one-to-one correspondence with the outer
corners of ${\mathbf T}_{\lambda}$. More precisely, if $(i, \lambda_i)$ is
an outer corner of ${\mathbf T}_{\lambda}$, then the product of the
polytopes $($simplices$)$ with vertices $\{x_1,\ldots,x_i\}$ and
$\{y_1,\ldots,y_{\lambda_i}\}$ is a facet of $X_{\lambda}$. Each
facet of $X_{\lambda}$ determines a rectangular region in the
Ferrers tableau ${\mathbf T}_{\lambda}$. The intersection of the
regions corresponding to two facets is again a rectangle that
corresponds to a product of smaller simplices. This product polytope
is the intersection of the two facets of $X_{\lambda}$. An example
is illustrated in Figure 4.
\begin{center}
\begin{pspicture}(-8,-1.5)(0,3.5)
\psset{xunit=.6cm, yunit=.6cm}

\pscustom[linestyle=solid,fillstyle=solid,fillcolor=lightgray]{
\psline(-7,3)(-5,3) \psline[liftpen=0](-5,1)(-7,1) }

\psline[linestyle=solid](-7,3)(-4,3)
\psline[linestyle=solid](-7,2)(-4,2)
\psline[linestyle=solid](-7,1)(-4,1)
\psline[linestyle=solid](-7,0)(-5,0)
\psline[linestyle=solid](-7,3)(-7,0)
\psline[linestyle=solid](-6,3)(-6,0)
\psline[linestyle=solid](-5,3)(-5,0)
\psline[linestyle=solid](-4,3)(-4,1)

\rput(-7.4,2.5){$x_1$} \rput(-7.4,1.5){$x_2$} \rput(-7.4,0.5){$x_3$}

\rput(-6.5,3.4){$y_1$} \rput(-5.5,3.4){$y_2$} \rput(-4.5,3.4){$y_3$}

\rput(-1.8,1){$\leadsto$}
\end{pspicture}
\begin{pspicture}(1,-2.5)(7,2.5)
\psset{xunit=.8cm, yunit=.8cm}

\pscustom[linestyle=dashed,fillstyle=solid,fillcolor=lightgray]{
\psline(2,2)(4.75,.5) \psline[liftpen=0](4.75,-1.5)(2,0) }

\pscustom[linestyle=solid,fillstyle=solid]{
\psline{*-*}(2,2)(4,2) \psline[liftpen=0](4,2)(4.75,.5) }

\psline[linestyle=dashed]{*-*}(2,0)(4,0)
\psline[linestyle=dashed](4,0)(4,2)

\psline[linestyle=solid]{*-*}(2,2)(2.75,-1.5)
\psline[linestyle=solid]{*-*}(2.75,-1.5)(2,0)
\psline[linestyle=solid]{*-*}(2,2)(2,0)

\psline[linestyle=dashed]{*-*}(2,2)(2,0)
\psline[linestyle=dashed](4,0)(5.5,-3)
\psline[linestyle=solid](2.75,-1.5)(5.5,-3)
\psline[linestyle=solid]{*-*}(4.75,0.5)(5.5,-3)

\psline[linestyle=solid](2,2)(2.75,-1.5)
\psline[linestyle=solid](2,2)(4.75,0.5)

\rput(4.75,-1.5){$\bullet$} \rput(4.75,.5){$\bullet$}

\rput(1.6,2.3){$x_1y_1$} \rput(1.6,-0.3){$x_2y_1$}
\rput(2.35,-1.85){$x_3y_1$}

\rput(3.5,2.3){$x_1y_3$} \rput(5.4,0.7){$x_1y_2$}
\rput(3.5,0.3){$x_2y_3$}

\rput(5.5,-3.3){$x_3y_2$} \rput(4.25,-1.75){$x_2y_2$}

\end{pspicture}
\end{center}
\begin{center}
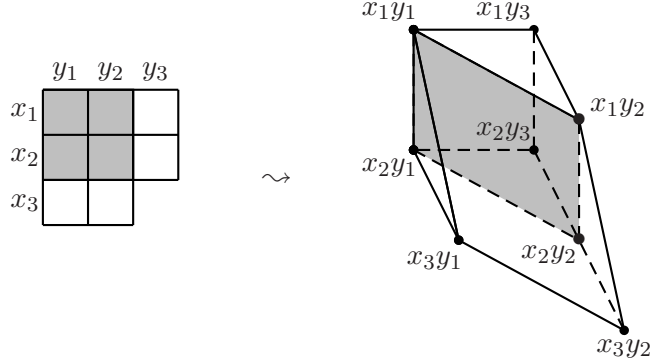
{\bf Figure 4: Faces of the polyhedral cell complex $X_{\lambda}$}
\end{center}

The main result of this section is:

\begin{Theorem}\label{exactness}
The complex ${\mathcal F}_{X_{\lambda}}$ provides the minimal free
$\ZZ^{m+n}$-graded resolution of $I_{\lambda}$.
\end{Theorem}

\noindent It is clear that ${\mathcal F}_{X_{\lambda}}$ also gives a
$\ZZ$-graded minimal free resolution of $I_{\lambda}$. In fact,
since the label of each $k$-dimensional face of $X_{\lambda}$ has
degree $k+2$ as noted above, we get $R^{F_k(X_{\lambda})} \cong
R^{f_k} (-k-2)$ where $f_k := |F_k(X_{\lambda})|$. Hence, we find
again $($as seen in Theorem~\ref{betti}$)$ that the $\ZZ$-graded
minimal free resolution of $I_{\lambda}$ is $2$-linear.

\medskip

\noindent{\bf Proof:} For fixed $m = \lambda_1$, we will induct on
$|\lambda| = \lambda_1 + \ldots + \lambda_n \geq m$. If $|\lambda| =
m$, then $X_{\lambda} = X_{1, m}$ and the claim follows from the
discussion above. Let $|\lambda| > m$. We divide the argument into
five steps:

{\bf (I)} For each $k \leq n+ \lambda_n -2$, let $G_{k-1} \subset
R^{F_k(X_{\lambda})}$ denote the free $R$-module generated by the
$k$-dimensional faces of $X_{\lambda}$ involving the vertex $x_n
y_{\lambda_n}$. Its rank is:
\[
{\rm rank} \, G_{k-1} = {n +\lambda_n-2 \choose k}.
\]
Indeed, each such face corresponds to the boxes lying on a suitable
grid of $a$ rows and $b$ columns indexed by $i_1 < i_2 < \cdots <
i_a=n$ and $j_1 < j_2 < \cdots < j_b = \lambda_n$, where $a+b-2=k$.
In this way, we see that the number of such $k$-dimensional faces
is:
\begin{eqnarray*}
{\rm rank} \, G_{k-1} & = & \sum_{a=1}^{k+1} {n-1 \choose
a-1}{\lambda_n-1
\choose b-1}  =  \sum_{a=1}^{k+1} {n-1 \choose a-1}{\lambda_n-1 \choose k-a+1} \\
& = & \sum_{j=0}^k {n-1 \choose j} {\lambda_n -1 \choose k-j} = {n
+\lambda_n-2 \choose k}.
\end{eqnarray*}
The latter equality in nothing but the Vandermonde convolution
\cite[Formula 3.1]{Gould}.

{\bf (II)} The proof of Theorem~\ref{betti} shows that there is a
partition $\lambda'$ such that $|\lambda'| = |\lambda| - 1$,
\[
I_{\lambda} = I_{\lambda'} + (x_ny_{\lambda_n}) \quad \mbox{and} \quad I_{\lambda'}
\colon x_ny_{\lambda_n} = (x_1, \ldots, x_{n-1}, y_1, \ldots,
y_{\lambda_n-1}).
\]
This provides the exact sequence:
\[
0 \rightarrow R/(x_1, \ldots, x_{n-1}, y_1, \ldots,
y_{\lambda_n-1})[-2] \stackrel{\cdot
x_ny_{\lambda_n}}{\longrightarrow} R/I_{\lambda'} \longrightarrow
R/I_{\lambda} \rightarrow 0.
\]

{\bf (III)} Let $\varepsilon$ be the incidence function of
$X_{\lambda}$ that gives the signs in ${\mathcal F}_{X_{\lambda}}$.
Then its restriction to $X_{\lambda'}$ is an incidence function too,
which we use to define the cell complex ${\mathcal
F}_{X_{\lambda'}}$.

Observe that, for each variable $l \in R$ and each non-empty face $P
\in X_{\lambda}$, there is a unique facet $Q$ of $P$ such that $m_P
= l \cdot m_Q$. We denote this facet $Q$ by $P/l$. Let $P$ denote an
$k$-dimensional face of $X_{\lambda}$ involving the monomial
$x_ny_{\lambda_n}$ and observe that $\partial_k(e_P)$ can be written
as:
\begin{eqnarray*}
\partial_k(e_P) & = & \sum_{l | m_P} l \cdot \varepsilon(P,P/l)e_{P/l} \\
& = & \sum_{\split{l \mid m_P}{l \nmid x_ny_{\lambda_n}}} l \cdot
\varepsilon(P, P/l) e_{P/l} + x_n \varepsilon(P,P/x_n) e_{P/x_n} +
y_{\lambda_n} \varepsilon(P,P/y_{\lambda_n}) e_{P/y_{\lambda_n}}
\\
& = & \varphi_{k-1}(e_P) + (-1)^k \delta_{k-1}(e_P),
\end{eqnarray*}
where
\begin{eqnarray*}
\varphi_{k-1}(e_P) & = & \sum_{\split{l \mid m_P}{l \nmid
x_ny_{\lambda_n}}} l \cdot \varepsilon(P, P/l)
e_{P/l} \\
\delta_{k-1}(e_P) & = & (-1)^k x_n \varepsilon(P,P/x_n) e_{P/x_n} +
(-1)^k y_{\lambda_n} \varepsilon(P,P/y_{\lambda_n})
e_{P/y_{\lambda_n}}.
\end{eqnarray*}
Note that $\varphi_{k-1}(e_P)$ is in $G_{k-2}$. Thus, we get a
sequence of graded $R$-modules:
\[
{\mathbb G}_{\bullet}: \qquad 0 \longrightarrow G_{n+\lambda_n - 3}
\stackrel{\varphi_{n+\lambda_n - 3}}{\longrightarrow} \ldots
\longrightarrow G_1 \stackrel{\varphi_{1}}{\longrightarrow} G_0
\longrightarrow 0,
\]
where the image of $\varphi_1$ is the ideal $(x_1, \ldots, x_{n-1},
y_1, \ldots, y_{\lambda_n-1})$. In fact, it is not too difficult to
see that ${\mathbb G}_{\bullet}$ is actually the Koszul complex on
$x_1, \ldots, x_{n-1}, y_1, \ldots, y_{\lambda_n-1}$ where the
degrees are shifted by $-2$.

{\bf (IV)} Set $H_k' = R^{F_k(X_{\lambda'})}$. Then
$R^{F_k(X_{\lambda})} = H_k' \oplus G_{k-1}$. Moreover, for each
generator $e_P \in G_{k-1}$, $\delta_{k-1}(e_P)$ is in $H_{k-1}'$.
Hence, we get the following square:
\[
\xymatrix{ G_{k-2}
\ar[r]^{\delta_{k-2}} & H'_{k-2} \\
G_{k-1} \ar[u]^{\varphi_{k-1}} \ar[r]_{\delta_{k-1}} & H'_{k-1}.
\ar[u]_{\partial'_{k-1}}}
\]
We claim that it is commutative, i.e.\ $\delta_{k-2} \circ
\varphi_{k-1} = \partial_{k-1}' \circ \delta_{k-1}$. Indeed, we have
that:
\begin{eqnarray*}
\lefteqn{\delta_{k-2}(\varphi_{k-1}(e_P))= } \\
& = & \delta_{k-2}\biggl( \sum_{l \mid \frac{m_P}{x_ny_{\lambda_n}}}
l \, \varepsilon(P, P/l) \,
e_{P/l} \biggr) \\
& = & (-1)^{k-1} \biggl( \sum_{l \mid \frac{m_P}{x_ny_{\lambda_n}}}
x_n l \, \varepsilon(P,P/l) \varepsilon(P/l, P/lx_n) e_{P/lx_n} +
y_{\lambda_n}l \, \varepsilon(P,P/l) \varepsilon(P/l,
P/ly_{\lambda_n}) e_{P/ly_{\lambda_n}} \biggr).
\end{eqnarray*}
On the other hand, we get:
\begin{eqnarray*}
\lefteqn{\partial_{k-1}'(\delta_{k-1}(e_P))= } \\
 & = & \partial_{k-1}' \biggl( (-1)^k
x_n \, \varepsilon(P,P/x_n)e_{P/x_n} + (-1)^k y_{\lambda_n} \,
\varepsilon(P, P/y_{\lambda_n}) e_{P/y_{\lambda_n}} \biggr) \\
& = & (-1)^k x_n \, \varepsilon(P, P/x_n) \sum_{l \mid
\frac{m_P}{x_n}} l \,
\varepsilon(P/x_n, P/lx_n) e_{P/lx_n}  \\
& & + (-1)^k y_{\lambda_n} \, \varepsilon(P, P/y_{\lambda_n})
\sum_{l \mid \frac{m_P}{y_{\lambda_n}}} l \,
\varepsilon(P/y_{\lambda_n},
P/ly_{\lambda_n}) e_{P/ly_{\lambda_n}} \\
& = & (-1)^k x_n \, \varepsilon(P, P/x_n) \biggl( \sum_{l \mid
\frac{m_P}{x_ny_{\lambda_n}}} l \, \varepsilon(P/x_n, P/lx_n)
e_{P/lx_n} + y_{\lambda_n} \, \varepsilon(P/x_n, P/x_ny_{\lambda_n})
e_{P/x_ny_{\lambda_n}} \biggr) \\
& & + (-1)^k y_{\lambda_n} \, \varepsilon(P, P/y_{\lambda_n})
\biggl( \sum_{l \mid \frac{m_P}{x_ny_{\lambda_n}}} l \,
\varepsilon(P/y_{\lambda_n}, P/ly_{\lambda_n}) e_{P/ly_{\lambda_n}}
+ x_n \, \varepsilon(P/y_{\lambda_n}, P/x_ny_{\lambda_n})
e_{P/x_ny_{\lambda_n}} \biggr) \\
& = & (-1)^k x_n \, \varepsilon(P, P/x_n) \biggl( \sum_{l \mid
\frac{m_P}{x_ny_{\lambda_n}}} l \, \varepsilon(P/x_n, P/lx_n)
e_{P/lx_n}
\biggr) \\
& & + (-1)^k y_{\lambda_n} \, \varepsilon(P, P/y_{\lambda_n})
\biggl( \sum_{l \mid \frac{m_P}{x_ny_{\lambda_n}}} l \,
\varepsilon(P/y_{\lambda_n}, P/ly_{\lambda_n})
e_{P/ly_{\lambda_n}} \biggr) \\
& = & (-1)^{k-1} \biggl( \sum_{l \mid \frac{m_P}{x_ny_{\lambda_n}}}
x_n l \,\varepsilon(P, P/l) \varepsilon(P/l, P/lx_n) e_{P/lx_n} +
y_{\lambda_n} l \, \varepsilon(P, P/l) \varepsilon(P/l,
P/ly_{\lambda_n}) e_{P/ly_{\lambda_n}} \biggr).
\end{eqnarray*}
Observe that the last two equalities follow from one of the properties of incidence
functions:
\[
\varepsilon(F,F/l) \cdot \varepsilon(F/l,F/lh)+
\varepsilon(F,F/h) \cdot \varepsilon(F/h,F/lh)=0.
\]

{\bf (V)} By Step (IV), $\delta_{\bullet}: {\mathbb G}_{\bullet} \to
{\mathcal F}_{X_{\lambda'}}$ is a morphisms of chain complexes. It
allows us to apply the mapping cone procedure to the exact sequence
in Step (II), which provides the desired free resolution of
$R/I_{\lambda}$. \QED

\medskip

In \cite{EGHP}, Eisenbud, Green, Hulek, and Popescu consider more
generally a projective subscheme that is the union of linear subspaces
and that has a $2$-linear free resolution. They construct a free
resolution of such a scheme $X$. However, it is not in general minimal
though it gives the exact number of minimal generators of the
homogeneous ideal $I_X$. Our Theorem~\ref{exactness} treats the
special case where $I_X$ is a monomial ideal, but our conclusion is
stronger.

\section{Characterization of Ferrers graphs}

In this brief section we establish an intrinsic characterization of
Ferrers graphs (not referring to a suitable labeling) by proving the
converse of Theorem~\ref{betti}. In other words, we characterize
Ferrers ideals as essentially the only edge ideals with a $2$-linear
free resolution among the ones arising from bipartite graphs. To
this end we will use a result of Fr\"oberg, which has been recently
refined in \cite{EGHP2}.

Let $G$ be a finite graph on the vertex set $V=\{ v_1, \ldots,
v_n\}$. We recall that the {\it complementary graph} $\overline{G}$
of $G$ is the graph $($on the same vertex set $V$ as $G)$ such that,
for vertices $v_i, v_j \in V$, the pair $(v_i, v_j)$ is an edge of
$\overline{G}$ if and only if $(v_i, v_j)$ is not an edge of $G$.
Furthermore, the graph $G$ is called {\it chordal} if every cycle of
$\overline{G}$ of length at least $4$ has a chord. With this
notation, a result of Fr\"oberg says:

\begin{Theorem}[Fr\"oberg~\mbox{\cite{Froberg}}]
The edge ideal of a graph $G$ has a $2$-linear free resolution if
and only if the complementary graph $\overline{G}$ is chordal.
\end{Theorem}

Note that adding an isolated vertex to a given graph $G$ does not
change the generating set nor the graded Betti numbers of the edge
ideal of $G$. Thus, it is harmless to assume that the graph does not
have isolated vertices. We are now ready to show:

\begin{Theorem}\label{converse}
Let $G$ be a bipartite graph without isolated vertices. Then its
edge ideal has a $2$-linear free resolution if and only if $G$ is
$($up to a relabeling of the vertices$)$ a Ferrers graph.
\end{Theorem}
\demo We have shown in Theorem~\ref{betti} that the condition is
sufficient. We now establish its necessity. Thus, let $G$ be a
bipartite graph on two distinct set of vertices, say $\{ x_1,
\ldots, x_n \}$ and $\{ y_1, \ldots, y_m \}$, and assume that its
edge ideal has a $2$-linear resolution. Let $\lambda_i$ be the
degree of $x_i$. By relabeling the vertices, we may also assume that
$m \geq \lambda_1 \geq \lambda_2 \geq \ldots \geq \lambda_n \geq 1$
and that the edges connected to $x_1$ are labeled $y_1, \ldots,
y_{\lambda_1}$.

For $1 \leq i \leq n$, we now claim that the $\lambda_i$ vertices
connected to $x_i$ are exactly the first consecutive $\lambda_i$
vertices contained in $\{ y_1, \ldots, y_{\lambda_{i-1}} \}$.
Indeed, there is nothing to prove if $i=1$. Let $i > 1$ and assume
that $x_i$ is connected to some $y_k$ with $k > \lambda_{i-1}$. Thus
there exists some  $j$, with $1 \leq j \leq \lambda_{i-1}$, such
that $x_iy_j$ is not an edge in $G$. Moreover, the induction
hypothesis provides that $x_{i-1}y_k$ is not an edge of $G$ either.
It follows that the complementary graph $\overline{G}$ contains the
cycle $\Gamma=\{ x_iy_j, x_{i-1}y_k, x_{i-1}x_i, y_jy_k \}$ of
length $4$. However, none of the chords of $\Gamma$, namely $x_iy_k$
and $x_{i-1}y_j$, belongs to $\overline{G}$. This contradicts
Fr\"oberg's theorem. Now, by relabeling the vertices of $G$ we may
assume that the vertices connected to $x_i$ are exactly the first
consecutive $\lambda_i$ vertices contained in $\{ y_1, \ldots,
y_{\lambda_{i-1}} \}$.

As a by-product of our argument, we observe that $\lambda_1$ is
exactly $m$. It also shows that if $(x_p, y_q)$ is an edge of $G$,
then so is $(x_h, y_k)$, provided $1 \leq h \leq p$ and $1 \leq k
\leq q$. In addition, $(x_1,y_m)$ and $(x_n,y_1)$ are edges of $G$.
Hence $G$ is a Ferrers graph as claimed. \QED

\medskip

The following example shows  that there are edge ideals with a
$2$-linear resolution which do not arise from a bipartite graph.
However, this ideal can be obtained as a specialization of a
suitable Ferrers ideal.

\begin{Example}\label{notbipartite}
{\rm Let $R=K[x_1, x_2, x_3]$ be a polynomial ring over a field $K$
and let $I$ be the edge ideal corresponding to a cycle of length
three, that is $I=(x_1x_2, \, x_1x_3, \, x_2x_3)$. Clearly, the
ideal $I$ does not arise from a bipartite graph. On the other hand,
it is an height two Cohen-Macaulay ideal with the following 2-linear
resolution
\[
0 \rightarrow R^{2}[-3] \stackrel{\varphi}{\longrightarrow}
R^{3}[-2] \longrightarrow R \longrightarrow R/I \rightarrow 0.
\]
We notice though that it can be obtained as a specialization of the
Ferrers ideal $I_{\lambda}=(x_1y_1, x_1y_2,$ $x_2y_1)$, by setting
$y_1:=x_3$ and $y_2:=x_2$ $($see \cite{CN2}$)$. }
\end{Example}

\medskip

Notice that the combination of Theorems~\ref{converse} and \ref{betti}
provides a complete description of the possible Betti numbers of edge
ideals of bipartite graphs with a $2$-linear free resolution.

\section{Toric rings associated to Ferrers ideals}\label{Toric}

Let ${\bf T} := {\bf T}_{\lambda}$ denote the Ferrers tableau
associated to a Ferrers
graph $G$ with partition $\lambda=(\lambda_1, \ldots, \lambda_s, 1,$
$\ldots, 1)$. We now define an associated tableau ${\bf T}'$
obtained from ${\bf T}$ by deleting all boxes in the first row
beyond the $\lambda_2$ one, and all boxed in the first column beyond
the $s$ one. Hence the partition $\lambda'$ associated to ${\bf T}'$
is $(\lambda_2, \lambda_2, \lambda_3, \ldots, \lambda_s)$. Observe
that, in this manner, the thickness of the outer border of ${\bf
T}'$ is at least $2$. From a combinatorial point of view, we removed
from $G$ all the vertices $($and, a fortiori, the corresponding
edges$)$ having degree $1$. An example is illustrated below:
\begin{center}
\begin{pspicture}(-5,-1)(10,4)
\psset{xunit=.6cm, yunit=.6cm}

\rput(0,-1){{\it Ferrers tableau ${\bf T}$}}

\rput(9,-1){{\it Ferrers tableau ${\bf T}'$}}

\rput(-2.5,4.5){$x_1$} \rput(-2.5,3.5){$x_2$} \rput(-2.5,2.5){$x_3$}
\rput(-2.5,1.5){$x_4$}\rput(-2.5,0.5){$x_5$}

\rput(-1.5,5.5){$y_1$} \rput(-0.5,5.5){$y_2$} \rput(0.5,5.5){$y_3$}
\rput(1.5,5.5){$y_4$}\rput(2.5,5.5){$y_5$}

\psline[linestyle=solid](-2,5)(3,5)
\psline[linestyle=solid](-2,4)(3,4)
\psline[linestyle=solid](-2,3)(2,3)
\psline[linestyle=solid](-2,2)(1,2)
\psline[linestyle=solid](-2,1)(0,1)
\psline[linestyle=solid](-2,0)(-1,0)

\psline[linestyle=solid](-2,5)(-2,0)
\psline[linestyle=solid](-1,5)(-1,0)
\psline[linestyle=solid](0,5)(0,1)
\psline[linestyle=solid](1,5)(1,2)
\psline[linestyle=solid](2,5)(2,3)
\psline[linestyle=solid](3,5)(3,4)

\rput(6.5,4.5){$x_1$} \rput(6.5,3.5){$x_2$} \rput(6.5,2.5){$x_3$}
\rput(6.5,1.5){$x_4$}

\rput(7.5,5.5){$y_1$} \rput(8.5,5.5){$y_2$} \rput(9.5,5.5){$y_3$}
\rput(10.5,5.5){$y_4$}

\psline[linestyle=solid](7,5)(11,5)
\psline[linestyle=solid](7,4)(11,4)
\psline[linestyle=solid](7,3)(11,3)
\psline[linestyle=solid](7,2)(10,2)
\psline[linestyle=solid](7,1)(9,1)

\psline[linestyle=solid](7,5)(7,1)
\psline[linestyle=solid](8,5)(8,1)
\psline[linestyle=solid](9,5)(9,1)
\psline[linestyle=solid](10,5)(10,2)
\psline[linestyle=solid](11,5)(11,3)
\psline[linestyle=solid](8,5)(8,1)

\rput(9,2){$\bullet$} \rput(10,3){$\bullet$}
\psline[linestyle=dashed](7.5,0.5)(11.5,4.5)

\psline[linewidth=1mm, linestyle=solid](-2,5)(-2,1)
\psline[linewidth=1mm, linestyle=solid](0,1)(0,2)
\psline[linewidth=1mm, linestyle=solid](1,2)(1,3)
\psline[linewidth=1mm, linestyle=solid](2,5)(2,3)

\psline[linewidth=1mm, linestyle=solid](-2.08,5)(2.08,5)
\psline[linewidth=1mm, linestyle=solid](-2.08,1)(0.08,1)
\psline[linewidth=1mm, linestyle=solid](-0.08,2)(1.08,2)
\psline[linewidth=1mm, linestyle=solid](0.92,3)(2.08,3)
\end{pspicture}
\end{center}
\begin{center}
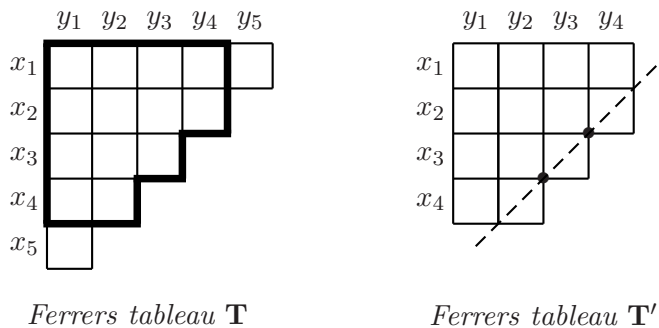
{\bf Figure 5: Ferrers tableaux ${\mathbf T}$ and ${\mathbf T}'$}
\end{center}
According to \cite{SVV}, the Rees algebra $R[It]$, the associated
graded ring ${\rm gr}_I(R)$ and the special fiber ring ${\mathcal
F}(I)$ of the edge ideal $I$ of every bipartite graph are normal
Cohen-Macaulay domains. Since edge ideals are generated in one
degree, the special fiber ring is also isomorphic to the toric ring
of the graph.

\begin{Proposition} \label{thm-ass-rings}
Let ${\bf X}=\{ x_1, \ldots, x_n \}$ and ${\bf Y} = \{ y_1, \ldots,
y_m \}$ be distinct sets of variables. Set $R=K[{\bf X}, {\bf
Y}]$, where $K$ is a field, and let $I_{\lambda}$ be the edge ideal
corresponding to a Ferrers graph $G_{\lambda}$ with associated
tableaux ${\bf T}$ and ${\bf T}'$, and partition
$\lambda=(\lambda_1, \lambda_2, \ldots, \lambda_s, 1, \ldots, 1)$.
Then the special fiber ring ${\mathcal F}(I_{\lambda})$ of
$I_{\lambda}$ has the following properties$:$
\begin{itemize}
\item[$(${\it a}$)$] ${\mathcal F}(I_{\lambda})$ is a Cohen-Macaulay normal
domain of dimension $n+m-1;$

\item[$(${\it b}$)$] ${\mathcal F}(I_{\lambda})$ is the ladder determinantal
ring $k[{\bf T}]/I_2({\bf T}');$

\item[$(${\it c}$)$] ${\mathcal F}(I_{\lambda})$ is Gorenstein if and only
if $\lambda_2=s$ and all the inside corners $($if any$)$ of the
Ferrers tableau ${\bf T}'$ lie on the main anti-diagonal of ${\bf
T}'$, i.e. $\{ (i,j) \in {\bf T}' \,|\, i+j=\lambda_2+1 \}$.
\end{itemize}
\end{Proposition}

\demo The result stated in $({\it a})$ is due to Simis,
Vasconcelos and Villarreal and holds for the special fiber ring of
the edge ideal of every  connected bipartite graph \cite{SVV}.
It is also recovered by part $({\it b})$, as ladder determinantal
rings are known to have such properties $($see \cite{N,HT,Cn}$)$.

In order to prove $(${\it b}$)$ observe that ${\mathcal
F}(I_{\lambda}) \cong K[x_iy_j\mbox{'s}]=K[G_{\lambda}]$, as
$I_{\lambda}$ is generated by homogeneous polynomials of the same
degree. Moreover, since $G_{\lambda}$ is a bipartite graph its
dimension is $m+n-1$ $($see \cite{SVV} or \cite[8.2.13]{vila}$)$.
Let ${\bf T}$ and ${\bf T}'$ denote the Ferrers tableaux associated
to $G_{\lambda}$. Let $T_{ij}$, for $(i,j)\in {\bf T}$, be
distinct variables: each variable is associated to the
corresponding box of the tableau ${\bf T}$ $($and ${\bf T}'$,
respectively$)$. By abuse of notation we also let ${\bf T}$ $($and ${\bf
T}'$, respectively$)$ denote the collection of these new variables.
We now consider the following epimorphism
\[
\pi\colon K[{\bf T}] \twoheadrightarrow K[G_{\lambda}] \cong
{\mathcal F}(I_{\lambda}),
\]
where $\pi(T_{ij})=x_iy_j$. We claim that the kernel of $\pi$ is
the determinantal ideal $I_2({\bf T}') = I_2({\bf T}') \cdot
k[{\bf T}]$ generated by the $2 \times 2$ minors of the one-sided
ladder ${\bf T}'$. It is clear that the ideal $I_2({\bf T}') \cdot
k[{\bf T}]$ is contained in the ideal ${\rm ker}(\pi)$. On the
other hand, we now show that these ideals have the same height.
Hence they coincide, as they are both prime ideals $($see \cite{N}
for the primeness of $I_2({\bf T}'))$. Indeed, we have
\begin{eqnarray*}
{\rm ht} \ {\rm ker}(\pi) & = & \dim k[{\bf T}] - \dim k[G_{\lambda}] \\
& = & (\lambda_1 + \lambda_2 + \cdots + \lambda_s + n-s) - (n+m-1) \\
& = & \lambda_2 + \cdots + \lambda_s -s +1
\end{eqnarray*}
$($as $\lambda_1=m)$, whereas
\begin{eqnarray*}
{\rm ht} \ I_2({\bf T}') \cdot k[{\bf T}] & = &
{\rm ht} \ I_2({\bf T}') \cdot k[{\bf T}'] \\
& = & \dim k[{\bf T}'] - \dim R_2({\bf T}') = \\
& = & \lambda_2 + \lambda_2 + \cdots + \lambda_s - (s+\lambda_2-1) \\
& = & \lambda_2 + \cdots + \lambda_s - s + 1.
\end{eqnarray*}
As far as $({\it c})$ is concerned, the Gorensteiness of ${\mathcal
F}(I_{\lambda})$ now follows from work of Conca \cite[2.5]{Cn}. \QED

\medskip

\begin{Corollary} \label{cor-toric-is-G}
The special fiber ring ${\mathcal F}(I_{\lambda})$ is Gorenstein if
and only if there is a partition $\mu$ such that the Ferrers ideal
$I_{\mu}$ is unmixed and there are variables such that the polynomial
rings over ${\mathcal F}(I_{\lambda})$ and  ${\mathcal
F}(I_{\mu})$, respectively, are isomorphic.
\end{Corollary}
\demo Assume that ${\mathcal F}(I_{\lambda})$ is Gorenstein. Then
define $\mu := (\lambda_2 + 1, \lambda_2,\ldots,\lambda_s,1) \in
\ZZ^{\lambda_2 + 1}$. It follows that the ladder determinantal ideals
determined by $\lambda$ and $\mu$, respectively, have the same
generators. Moreover, Proposition~\ref{thm-ass-rings} and
Corollary~\ref{unmixed} provide that $I_{\mu}$ is unmixed.

Conversely, if $I_{\mu}$ is unmixed, then we see that ${\mathcal
  F}(I_{\mu})$ is Gorenstein.
 \QED

\bigskip

As announced earlier, we now turn our attention to the computation
of the Hilbert series of the toric ring $K[G_{\lambda}]$:
This is a
highly investigated area of research, see, for example, \cite{Abhy, BH, Cn,
  Cn2, CnH, G,
HT, K, KM, KP, KR, Kulkarni, N, R, wang}. While most of these works
involve --- to a different extent --- path counting arguments, we
offer here a new and self-contained approach based on Gorenstein
liaison theory.
Proposition~\ref{thm-ass-rings} implies that, for each partition
$\lambda \in \NN^n$, there is a unique polynomial $p_{\lambda} \in
\ZZ[t]$ such that the Hilbert series of $K[G_{\lambda}] \cong
{\mathcal F} (I_{\lambda})$ can be written as:
\[
P(K[G_{\lambda}], t) = \frac{p_{\lambda} (t)}{(1-t)^{n+m-1}}.
\]
Note that the multiplicity of $K[G_{\lambda}]$ is $e(K[G_{\lambda}])
= p_{\lambda}(1)$. With this notation and using Gorenstein liaison theory
methods, we establish the following key result, which provides a
simple recursive formula for the Hilbert series.

\begin{Lemma} \label{lemma-liai-rec}
Let $\lambda = (\lambda_1,\ldots,\lambda_n) \in \NN^n$ be a
partition with $\lambda_n \geq 2$. Set $\lambda'' =
(\lambda_1,\ldots,\lambda_{n-1}, \lambda_n - 1) \in \NN^n$ and
$\lambda' = (\lambda_1 - \lambda_n +1,\ldots,\lambda_{n-1} -
\lambda_n + 1) \in \NN^{n-1}$. If $n \geq 3$, then there is the
following relation among Hilbert series$:$
\[
p_{\lambda} (t) = p_{\lambda''} (t)  + t \cdot p_{\lambda'} (t).
\]
\end{Lemma}

\demo We need some more notation. Given the partition $\lambda \in
\ZZ^n$, we define $S := K[{\bf T}]$ as the polynomial ring in the
$\lambda_1 + \ldots + \lambda_n$ variables $T_{i j}$ and the ideal
$J_{\lambda} \subset S$ by $S/J_{\lambda} := K[G_{\lambda}]$. Let ${\bf T''}$
be the Ferrers tableau associated to $\lambda''$ and let ${\bf
\widetilde{T}}$ be the Ferrers tableau to the partition
$\widetilde{\lambda} := (\lambda_1,\ldots,\lambda_{n-1}) \in
\ZZ^{n-1}$. Furthermore, denote by $N$ and ${\bf
\widetilde{\widetilde{T}}}$ the subtableaux of ${\bf \widetilde{T}}$
consisting of the first $(\lambda_n - 1)$ and the remaining columns,
respectively. Let $I_1 (N)$ be the ideal generated by the entries of
$N$ and let $I_2 ({\bf \widetilde{\widetilde{T}}})$ be the ideal
generated by the $2 \times 2$ minors whose entries are in ${\bf
\widetilde{\widetilde{T}}}$. Finally, let $V, W, V' \subset \Proj\,
(S)$ be the subvarieties that are defined by $J_{\lambda},
J_{\lambda''}$, and $J_{\widetilde{\lambda}} + I_1 (N)$,
respectively.

In \cite[proof of Theorem 2.1]{Gorla}, Gorla shows that $V$ is an
elementary biliaison of $V'$ on $W$. Thus, $V$ is linearly
equivalent to the basic double link of $V'$ on $W$. In particular,
both have the same Hilbert function. If follows (see, for instance,
\cite[Lemma 4.8]{KMMNP}) that the Hilbert functions satisfy for all
integers $j$:
\[
h_V (j) = h_{V'} (j-1) + h_W (j) - h_W (j-1).
\]
In terms of Hilbert series this reads as:
\[
P(S/J_{\lambda}, t) = t \cdot P (S/I_{V'}, t)   + (1-t) \cdot P( S/J_{\lambda''} S, t).
\]
Thus we get using Proposition \ref{thm-ass-rings}:
\begin{equation} \label{eq-h-ser}
\frac{p_{\lambda} (t)}{(1-t)^{m+n-1}} = t \cdot P (S/I_{V'}, t) +
(1-t) \cdot \frac{p_{\lambda''} (t)}{(1-t)^{m+n}}
\end{equation}
because ${\displaystyle P( S/J_{\lambda''} S, t) =  \frac{1}{1-t}
\cdot P(K[G_{\lambda''}], t) = \frac{p_{\lambda''}
(t)}{(1-t)^{m+n}}}$.

The definition of the homogeneous ideal of $V'$ implies:
\[
I_{V'} = I_2 ({\bf \widetilde{\widetilde{T}}}) + I_1 (N).
\]
It follows that $S/I_{V'}$ is isomorphic to a polynomial ring in
$\lambda_n$ variables over  $K[{\bf \widetilde{\widetilde{T}}}]/I_2
({\bf \widetilde{\widetilde{T}}}) \cong K[G_{\lambda'}]$. Since
$\dim K[G_{\lambda'}] = m+n - \lambda_n - 1$, we get:
\[
P (S/I_{V'}, t) = \frac{1}{(1-t)^{\lambda_n}} \cdot P(K[G_{\lambda'}], t)
= \frac{1}{(1-t)^{\lambda_n}} \cdot \frac{p_{\lambda'} (t)}{(1-t)^{m+n-\lambda_n - 1}}
=  \frac{p_{\lambda'} (t)}{(1-t)^{m+n- 1}}.
\]
Substituting in Equation (\ref{eq-h-ser}), the claim follows. \QED

\bigskip

As a first consequence, we derive an explicit formula for the Hilbert
series. Observe that
all terms are non-negative.

\begin{Theorem} \label{thm-Hilb-series}
Let $\lambda = (\lambda_1,\ldots,\lambda_n) \in \NN^n$ be a partition
with $n \geq 2$.
Then the numerator of the normalized Hilbert series of
$K[G_{\lambda}]$ is:
\[
p_{\lambda} (t) = 1 + h_1 (\lambda) \cdot t + \cdots + h_{n-1}
(\lambda) \cdot  t^{n-1},
\]
where
\begin{equation} \label{eq-h1}
h_1 (\lambda) = \sum_{j=2}^n (\lambda_j - 1)
\end{equation}
and
\begin{equation} \label{eq-hk}
h_k = (\lambda) \sum_{2 \leq i_1 < i_2 < \ldots < i_k \leq n} \
\sum_{j_{k-1} =  \lambda_{i_1} - \lambda_{i_k} - k +
2}^{\lambda_{i_1} - k} \ \sum_{j_{k-2} =  \lambda_{i_1} -
\lambda_{i_{k-1}} - k + 3}^{j_{k-1}} \ldots \sum_{j_1 =
\lambda_{i_1} - \lambda_{i_{2}}}^{j_{2}} j_1,
\end{equation}
for $k \geq 2$.
\end{Theorem}

\demo We use the notation introduced in Lemma \ref{lemma-liai-rec}
and its proof. This result implies for all $k \in \NN$:
\begin{equation} \label{eq-recur}
h_k (\lambda) = h_k (\lambda'') + h_{k-1} (\lambda').
\end{equation}
It is easy to see that this recursion provides the formula for $h_1 (\lambda)$.

We know induct on $n \geq 2$. If $n = 2$, the minimal free
resolution of ${\mathcal F}(I_{\lambda})$ is given by an
Eagon-Northcott complex. This implies in particular that
\[
p_{\lambda} (t) = 1 + (\lambda_2 - 1) \cdot t,
\]
as claimed. Let $n \geq 3$. Now we induct on $k \geq 2$. Since the
case $k = 2$ is similar, but easier than the general case, we assume
$k \geq 3$. We now induct on $\lambda_n \geq 1$. If $\lambda_n = 1$,
then the sum ${\displaystyle \sum_{j_{k-1} = \lambda_{i_1} -
\lambda_n - k + 2}^{\lambda_{i_1} - k}}$ vanishes. Thus in the
formula for $h_k (\lambda)$ all sums with $i_k = n$ vanish. This
implies that we have to show $h_k (\lambda) = h_k
(\widetilde{\lambda})$ where $\widetilde{\lambda} =
(\lambda_1,\ldots,\lambda_{n-1})$. But this is true because the
ideals $J_{\lambda}$ and $J_{\widetilde{\lambda}}$ have the same
generators.

Finally, we may assume that $\lambda_n \geq 2$. Then the induction
hypotheses and Formula (\ref{eq-recur}) provide by distinguishing
the cases $i_k < n$ and $i_k = n$:
\begin{eqnarray*}
h_k (\lambda) & = & \sum_{2 \leq i_1 < \ldots < i_{k} \leq n-1} \
\sum_{j_{k-1} =  \lambda_{i_1} - \lambda_{i_k} - k +
2}^{\lambda_{i_1} - k} \ \sum_{j_{k-2} =  \lambda_{i_1} -
\lambda_{i_{k-1}} - k + 3}^{j_{k-1}} \ldots
\sum_{j_1 =  \lambda_{i_1} - \lambda_{i_{2}}}^{j_{2}} j_1 \\
& & + \sum_{2 \leq i_1 < \ldots < i_{k-1} \leq n-1} \ \sum_{j_{k-1}
= \lambda_{i_1} - \lambda_{{n}} - k + 3}^{\lambda_{i_1} - k} \
\sum_{j_{k-2} =  \lambda_{i_1} - \lambda_{i_{k-1}} - k +
3}^{j_{k-1}} \ldots
\sum_{j_1 =  \lambda_{i_1} - \lambda_{i_{2}}}^{j_{2}} j_1 \\
& & + \sum_{2 \leq i_1 < \ldots < i_{k-1} \leq n-1} \ \sum_{j_{k-2}
= \lambda_{i_1} - \lambda_{i_{k-1}} - k + 3}^{\lambda_{i_1} -
\lambda_n - k + 2} \ \sum_{j_{k-3} =  \lambda_{i_1} -
\lambda_{i_{k-2}} - k + 4}^{j_{k-2}} \ldots
\sum_{j_1 =  \lambda_{i_1} - \lambda_{i_{2}}}^{j_{2}} j_1 \\
& = & \sum_{2 \leq i_1 < \ldots < i_{k} \leq n-1} \ \sum_{j_{k-1} =
\lambda_{i_1} - \lambda_{i_k} - k + 2}^{\lambda_{i_1} - k} \
\sum_{j_{k-2} =  \lambda_{i_1} - \lambda_{i_{k-1}} - k +
3}^{j_{k-1}} \ldots
\sum_{j_1 =  \lambda_{i_1} - \lambda_{i_{2}}}^{j_{2}} j_1 \\
& & + \sum_{2 \leq i_1 < \ldots < i_{k-1} \leq n-1} \ \sum_{j_{k-1}
= \lambda_{i_1} - \lambda_{{n}} - k + 2}^{\lambda_{i_1} - k} \
\sum_{j_{k-2} =  \lambda_{i_1} - \lambda_{i_{k-1}} - k +
3}^{j_{k-1}} \ldots
\sum_{j_1 =  \lambda_{i_1} - \lambda_{i_{2}}}^{j_{2}} j_1 \\
& = & \sum_{2 \leq i_1 < i_2 \ldots < i_k \leq n} \ \sum_{j_{k-1} =
\lambda_{i_1} - \lambda_{i_k} - k + 2}^{\lambda_{i_1} - k} \
\sum_{j_{k-2} =  \lambda_{i_1} - \lambda_{i_{k-1}} - k +
3}^{j_{k-1}} \ldots \sum_{j_1 =  \lambda_{i_1} -
\lambda_{i_{2}}}^{j_{2}} j_1.
\end{eqnarray*}
This completes the proof. \QED

\begin{Remark} \label{rem-counting} { \rm It is well-known that the
    coefficient $h_i (\lambda)$ of $t^i$ in $p_{\lambda} (t)$ has a
    combinatorial interpretation. In fact, interpreting the Ferrers
    tableau as a bounded region in the lattice $\ZZ^2$, $h_i
    (\lambda)$ is the number of lattice paths inside ${\bf
    T}_{\lambda}$ that start in the south-west corner, end in the
    north-east corner, and have exactly $i$ east-north turns (see
    \cite{Abhy, GV, HT, KP}).
}
\end{Remark}

\medskip

For the multiplicity of $K[G_{\lambda}]$ we obtain a somewhat simpler formula:

\begin{Corollary} \label{cor-mult}
\[
e(K[G_{\lambda}]) = \sum_{j_{n-2} =  \lambda_{2} - \lambda_{n}
+1}^{\lambda_{2}} \sum_{j_{n-3} =  \lambda_{2} - \lambda_{n-1} +
1}^{j_{n-2}} \ldots \sum_{j_1 =  \lambda_{2} - \lambda_{3} +
1}^{j_{2}} j_1.
\]
\end{Corollary}

\demo Since $e({\mathcal F}(I_\lambda)) = p_{\lambda} (1)$, this
follows from Theorem \ref{thm-Hilb-series}. However,
computationally, it is easier to use more directly Lemma
\ref{lemma-liai-rec} which implies $e(K[G_{\lambda}]) =
e(K[G_{\lambda'}]) + e(K[G_{\lambda''}])$. \QED

\medskip

The method of proof, using Gorenstein liaison theory, applies to {\em
  all} ladder determinantal ideals  \cite{Gorla}. However, here we
  restrict ourselves to the ideals related to Ferrers graphs, i.e.\ to
  one-sided ladder determinatal ideals generated by $2 \times 2$
  minors.
\smallskip

We recall that for a finitely generated graded module $M$ (over an
affine $K$-algebra) a suitable measure for the complexity of its
resolution (hence of $M$ itself) is given by the Castelnuovo-Mumford
regularity ${\rm reg}(M)$, that is $\max\{j-i \s \beta_{ij} \not=0
\}$, where $\beta_{ij}$ are the graded Betti numbers of $M$. On the
other hand the $a$-invariant $a(M)$ of $M$ is the degree of the Hilbert
series of $M$ as a rational function. In general these numbers are
related by $a(M) \leq {\rm reg}(M) - {\rm depth}(M)$, with equality
if $M$ is Cohen-Macaulay. In the latter case, we thus have that
$\reg (M)$ equals the degree of the numerator $p_M(t)$ of the
Hillbert series of $M$. One can also interpret $a(M)$ in terms of
non-vanishing of the top local cohomology of $M$. The approach we
followed thus far allows us to easily compute the
Castelnuovo-Mumford regularity and the $a$-invariant of the toric
ring $K[G_{\lambda}]$. Before stating our results, we recall that
for a partition $\lambda$ we set $s := s (\lambda) := \lambda_2^*$.
Note that $\lambda_s \geq 2$.

\begin{Proposition} \label{prop-reg}
Let $\lambda$ be a partition such that $\lambda_2 \geq 2$. Then the
Castelnuovo-Mumford regularity of the toric ring of the Ferrers
graph $G_{\lambda}$ is:
\begin{eqnarray*}
\reg (K[G_{\lambda}]) & = &  \min \{\lambda_2^* - 1, \{\lambda_j + j
- 3 \s 2 \leq j \leq \lambda_2^* =:s \}\} \\ & = & \left
\{\begin{array}{ll}
s - 1 & \mif \lambda_s \geq 3 \\
\min \{ j-1 \s \lambda_j = 2 \} & \mif \lambda_s = 2.
\end{array} \right.
\end{eqnarray*}
\end{Proposition}

\demo The second equality follows simply by evaluating the minimum
using $\lambda_s + s - 3 > s - 1$ if $\lambda_s \geq 3$. In order to
show the first equality, we note that the Cohen-Macaulayness of
$K[G_{\lambda}]$ implies that $\reg (K[G_{\lambda}])= \deg
p_{\lambda}(t)$. Now for this proof denote by $r_{\lambda} - 1$ the
right-hand side of the claim. Then we have to show that $\deg
p_{\lambda} = r_{\lambda} - 1$. This follows directly from
Theorem~\ref{thm-Hilb-series}. Alternatively, we can use
Lemma~\ref{lemma-liai-rec}, and it suffices to show (using its
notation):
\[
r_{\lambda} = \max \{ r_{\lambda''}, 1 + r_{\lambda'} \}.
\]
But this can be easily checked. \QED

\medskip

\begin{Corollary} \label{invariants}
Let $\lambda$ be a partition such that $\lambda_2 \geq 2$. Then the
$a$-invariant of the toric ring of the Ferrers graph $G_{\lambda}$
is$:$
\[
a(K[G_{\lambda}])= -(n+m-1) + \min \{\lambda_2^* - 1, \{\lambda_j +
j - 3 \s 2 \leq j \leq \lambda_2^* \}\}
\]
\end{Corollary}
\demo Our claim follows from Theorem~\ref{thm-ass-rings}, the
equality $a(K[G_{\lambda}])= - {\rm dim}(K[G_{\lambda}]) + {\rm
reg}(K[G_{\lambda}])$ and Proposition~\ref{prop-reg}. \QED

\begin{Remark}{\rm
We find it noteworthy, at this stage, to highlight an existing
connection between $a$-invariants and Integer Programming
techniques, as pointed out by Valencia and Villarreal in
\cite{ValVila}. In fact, for the type of ideals considered in this
paper, the computation of the $a$-invariant amounts to the
computation of the maximum number of edge disjoint directed cuts or
equivalently to the minimum cardinality of the edge set that
contains at least one edge of each directed cut. The latter is not
easily computable using combinatorial techniques, thus it is
remarkable that
 Corollary~\ref{invariants}  provides  an
explicit formula. }
\end{Remark}

We conclude this section by discussing  particular classes of Ferrers
graphs where the formulas simplify considerably. In
Example~\ref{genericnbym} we recover in a simple way the expression
for the multiplicity (due to Herzog and Trung \cite{HT}) and the
coefficients of the Hilbert series of the $2 \times 2$ minors of a
generic $n \times m$ matrix (due to Conca and Herzog \cite{CnH}). In
the same simple fashion, we recover in Example~\ref{cor-class} a
result that appears in \cite{wang2}.

\begin{Example} \label{genericnbym} {\rm
Let $2 \leq n , m $ be integers and consider the partition $\lambda
= (\lambda_1,\ldots,\lambda_n)$ where $ \lambda_i := m$, i.e.\
$G_{\lambda}$ is the complete bipartite graph ${\mathcal K}_{n,m}$. Then the
coefficients of the polynomial in the Hilbert series of the toric
ring $K[G_{\lambda}]$ as well as its multiplicity are$:$
\[
h_k(\lambda) = {m-1 \choose k}{n-1 \choose k} \qquad \mbox{and} \qquad
e(K[G_{\lambda}]) = {n+m-2 \choose m-1}.
\]}
\end{Example}
\demo In the calculations that will follow we will make a repeated
use of the combinatorial identity
\[
\sum_{j_t=1}^{j_{t+1}}{ j_t+t-1 \choose j_t-1} = { j_{t+1} + (t+1)-1
\choose j_{t+1}-1},
\]
which can be found in \cite[Formula 1.49]{Gould}. According to
Theorem~\ref{thm-Hilb-series} we only need to show  the
formula for $h_k$ when $k \geq 2$, as in the other two cases the
expression is trivially verified. In this particular case the
expression in Theorem~\ref{thm-Hilb-series} reduces to:
\begin{eqnarray*}
h_k (\lambda) & = & \sum_{2 \leq i_1 < i_2 < \ldots < i_k \leq n} \
\sum_{j_{k-1} = - k + 2}^{m - k} \ \sum_{j_{k-2} = - k +
3}^{j_{k-1}} \ldots \sum_{j_2=-1}^{j_3} \ \sum_{j_1 = 0}^{j_{2}} j_1
\\
& = & \sum_{2 \leq i_1 < i_2 < \ldots < i_k \leq n} \ \sum_{j_{k-1}
= -k+2}^{m - k} \ \sum_{j_{k-2} = -k+3}^{j_{k-1}} \ldots
\sum_{j_2=-1}^{j_3} {j_2+1 \choose 2}  \\
& = & \sum_{2 \leq i_1 < i_2 < \ldots < i_k \leq n} \ \sum_{j_{k-1}
= -k+2}^{m - k} \ \sum_{j_{k-2} = -k+3}^{j_{k-1}} \ldots
\sum_{j_3=-2}^{j_4} {j_3+2 \choose 3}  \\
& = & \sum_{2 \leq i_1 < i_2 < \ldots < i_k \leq n} \ \sum_{j_{k-1}
= -k+2}^{m - k} {j_{k-1} + (k-1) -1 \choose j_{k-1}-1} \\
& = & \sum_{2 \leq i_1 < i_2 < \ldots < i_k \leq n} {m-1 \choose k}
 =  {m-1 \choose k}{n-1 \choose k}.
\\
\end{eqnarray*}
Finally, the expression for the multiplicity of $K[G_{\lambda}]$
follows from Corollary~\ref{cor-mult} by performing similar
computations. \QED

\medskip

\begin{Example} \label{cor-class} {\rm
Let $2 \leq n \leq m $ be integers and consider the partition
$\lambda = (\lambda_1,\ldots,\lambda_n)$ where $ \lambda_k :=
m+1-k$. Then the Hilbert series of the toric ring $K[G_{\lambda}]$
is$:$
\[
P(K[G_{\lambda}], t) = \frac{1 + h_1 (\lambda) \cdot t + \cdots +
h_{n-1} (\lambda) \cdot  t^{n-1}}{(1-t)^{m+n-1}},
\]
where
\[
h_k (\lambda) = \binom{n-1}{k} \binom{m-2}{k} - \binom{n-1}{k+1}
\binom{m-2}{k-1}
\]
for all $k = 0, \ldots, n-1$. In particular, the multiplicity is$:$
\[
e(K[G_{\lambda}]) = \frac{m-n+1}{m} \cdot \binom{m+n-2}{m-1}.
\] }
\end{Example}
\demo We induct on $n \geq 2$. Theorem~\ref{thm-Hilb-series}
immediately provides the claim about $h_0 (\lambda)$ and $h_1
(\lambda)$. Thus we may assume $n \geq 3$. Consider the partition
$\bar{\lambda} = (m,\ldots,m+2-n, j)$ where $1 \leq j \leq m+1-n$,
i.e., $\bar{\lambda}$ differs from the given partition $\lambda$ at
most in the last entry. Then, using the notation of Lemma
\ref{lemma-liai-rec}, if $j \geq 2$,  we get $\bar{\lambda''} =
(m,\ldots,m+2-n, j-1) \in \NN^{n}$ and $\bar{\lambda'} = (m - j
+1,\ldots,m-j+3-n) \in \NN^{n-1}$. Note that the induction
hypothesis applies to $\bar{\lambda'}$. Letting $j$ vary between $1$
and $m+1-n$, Lemma \ref{lemma-liai-rec} and the induction hypothesis
provide:
\begin{eqnarray*}
h_k (\lambda) & = & h_k (m,m-1,\ldots,m+2-n) + \sum_{j=2}^{m+1-n}
h_{k-1} (m-j+1,\ldots,m-j+3-n) \\
& = & \binom{n-2}{k} \binom{m-2}{k} - \binom{n-2}{k+1} \binom{m-2}{k-1} + \\
& & \sum_{j=2}^{m+1-n} \left [\binom{n-2}{k-1} \binom{m-j-1}{k-1} -
\binom{n-2}{k} \binom{m-j-1}{k-2} \right ] \\[1ex]
& = & \binom{n-2}{k} \binom{m-2}{k} - \binom{n-2}{k+1} \binom{m-2}{k-1} + \\
& & \binom{n-2}{k-1} \cdot \left [ \binom{m-2}{k} - \binom{n-2}{k}
\right ] -
\binom{n-2}{k} \cdot \left [ \binom{m-2}{k-1} - \binom{n-2}{k-1} \right ] \\[1ex]
& = & \binom{n-1}{k} \binom{m-2}{k} - \binom{n-1}{k+1}
\binom{m-2}{k-1},
\end{eqnarray*}
as claimed. Finally, as noted earlier, $e(K[G_{\lambda}]) =
p_{\lambda}(1)$. Thus, using \cite[Formula 3.20]{Gould}, we get
\begin{eqnarray*}
e(K[G_{\lambda}]) & = & \sum_{k=0}^{n-1} \left [ \binom{n-1}{k}
\binom{m-2}{k} - \binom{n-1}{k+1} \binom{m-2}{k-1} \right ]
\\[1ex]
& = & \binom{m+n-3}{n-1} - \binom{m+n-3}{n-3} \\[1ex]
& = & \frac{m-n+1}{m} \cdot \binom{m+n-2}{m-1},
\end{eqnarray*}
where the last equality is easy to verify. \QED

\medskip

As reflected in the coefficients $h_k$'s, one should observe that
the roles of $m$ and $n$ are not symmetric in the previous
corollary, since we consider the partition $\lambda=(m, m-1,
\ldots, m-n+1)$. However, in the case $m = n$, the above
formul\ae\/ greatly simplify (becoming symmetric!) and the rings
have particularly good properties:

\begin{Example}{\rm
Consider the partition $\lambda = (n, n-1,\ldots,2, 1) \in \ZZ^n$.
The ring $R/I_{\lambda}$ cogenerated by the edge ideal $I_{\lambda}$
of the associated Ferrers graph $G_{\lambda}$ is Cohen-Macaulay (by
{\rm Corollary~\ref{CM}}) with multiplicity $n+1$ and Hilbert
series:
\[
P(R/I_{\lambda},t) = \frac{1+nt}{(1-t)^n}.
\]
The toric ring $K[G_{\lambda}]$ is Gorenstein with Hilbert series:
\[
P(K[G_{\lambda}], t) = \frac{\displaystyle\sum_{k=0}^{n-2}
\frac{\displaystyle{n-2 \choose k } {n-1 \choose k} }{k+1} \
t^k}{(1-t)^{2n-1}}.
\]
In particular, the  multiplicity of $K[G_{\lambda}]$ is the {\it
Catalan number}:
\[
e(K[G_{\lambda}]) = \displaystyle\frac{\displaystyle{2(n-1) \choose
n-1}}{n}.
\]
We refer the interested reader to \cite{St, St2} for a wealth of
information about Catalan numbers.}
\end{Example}

\end{document}